\documentclass{amsart}

\usepackage{geometry} 
\usepackage{stmaryrd}
\usepackage{graphicx}
\usepackage{amssymb}
\usepackage{epstopdf}
\usepackage{booktabs}
\usepackage{enumerate}
\usepackage{multirow}
\usepackage{color}

\newcommand{\bld}[1]{\boldsymbol{#1}}
\newcommand{\uhat}{{\widehat{u}}}
\newcommand{\qhat}{{\widehat{{\bld{q}}}}}
\newcommand{\what}{{\widehat{w}}}
\newcommand{\dO}{{\partial\Omega}}
\newcommand{\Oh}{{\mathcal{T}_h}}
\newcommand{\dOh}{{\partial \Oh}}

\definecolor{red}{rgb}{1,0,0}

\begin{document}

\title[The WG methods are rewritings of the HDG methods]{The Weak Galerkin methods are rewritings of the Hybridizable Discontinuous Galerkin methods}

\author{Bernardo Cockburn}
\address{School of Mathematics, 127 Vincent Hall, University of Minnesota, Minneapolis,
Minnesota 55455}
\email{cockburn@math.umn.edu} 
\date{\today}                                           

\maketitle

\begin{abstract}
We establish that the Weak Galerkin methods are rewritings of the hybridizable discontinuous Galerkin methods.%
\end{abstract}

\section{Introduction}
Nowhere in the by now massive literature of Weak Galerkin (WG) methods, there has been a serious effort to 
clarify the relation between the WG and the Hybridizable Discontinuous Galerkin (HDG) methods. The need for  clarification stems
from the fact that the two methods are defined by using different weak formulations.
To the outsider, this has produced the impression that the methods are actually different.
In this paper, we dispel such an impression by showing that the WG methods are nothing but rewritings of the HDG methods. 

The paper is organized as follows. In Section 2, we begin by establishing this fact in the framework of steady-state diffusion, thus expanding previous work presented in Subsection \S 6.6 in \cite{CockburnDurham16}. In Section 3, we briefly indicate the equivalence of the methods in the framework of Stokes flow. In Section 4, we consider, in full detail, the case of the bilaplacian.
In Section 5, we critically assess several questionable statements made in the literature of the WG methods. We end with some concluding remarks.

\section{Second-order elliptic problems}
Here, we establish the equivalence between the HDG and the WG methods in the framework of the 
model elliptic problem
\begin{alignat*}{2}
\mathrm{c}\;{\bld q}=&-\nabla u&&\qquad\mbox{in }\Omega,
\\
\nabla\cdot {\bld q} =&\;f&&\qquad\mbox{in }\Omega,
\\
u=&\;u_{D}&&\qquad\mbox{on }\partial\Omega.
\end{alignat*}
We consider {\em four} WG methods for this model problem. We have to do this 
because the WG methods have been characterized in different ways in different papers. 

We begin by establishing that the first WG methods, proposed in 2013 \cite{WangYeMain13}, 
are a rewriting of the hybridized version 
of mixed methods like the 1994 mixed method introduced in \cite{Chen94} and the 1997 mixed methods in 
\cite{ArbogastWheelerYotov97} which were also called {\em expanded mixed methods} therein. The mixed methods fit the general formulation of the HDG methods and the only difference between them is that, unlike the HDG methods, the mixed methods do not need stabilization.

We  then consider the 
WG method proposed in 2014 in \cite{WangYemixedsecondorder14}, 
the WG method proposed in 2015 in \cite{MuWangYepolytopal15}, 
and the WG method proposed in 2015 in \cite{MuWangYeLS15}. 
We show that each of these three WG methods are simple rewritings of HDG methods.

\subsection{The HDG methods}
\subsubsection{The original formulation of the HDG methods}
Let us recall the definition of the
HDG methods as introduced for the first time in 2009 in \cite{CockburnGopalakrishnanLazarov09}.
The (conforming) mesh
of the domain $\Omega$ will be denoted by $\Oh$. The set of all the faces $F$ of the
elements of the mesh is denoted by $\mathcal{F}_h$. The set
$\partial\Oh:=\{\partial K:\; K\in \Oh\}$ is the set of the boundaries of all
the elements $K$ of the triangulation $\Oh$.

We also write 
\begin{alignat*}{2}
&(\mathrm{c}\,{{\bld q}_h},{\bld v})_{\Oh}:= \sum_{K\in\Oh}(\mathrm{c}\,{{\bld q}_h},{\bld
  v})_{K},
&&\quad\mbox{ where }(\mathrm{c}\,{{\bld q}_h},{\bld v})_K:=\int_K\,\mathrm{c}\, {{\bld
    q}_h}\cdot{\bld v},
\\
&\langle {\mu},{\bld v} \cdot {\bld n}\rangle_{\dOh}
:=\sum_{K\in\Oh}\langle {\mu},{\bld v} \cdot {\bld n}\rangle_{\partial K}
&&\quad\mbox{ where }
\langle {\mu},{\bld v} \cdot {\bld n}\rangle_{\partial K}
:=\int_{\partial K} {\mu}\, {\bld v} \cdot {\bld n}.
\end{alignat*}
Similar notation holds for other integrals. 

The HDG method 
seeks an approximation $({\bld{q}_h}, {{u_h}},{\uhat_h})$,   to the exact solution 
\[
({\bld{q}}:=-\mathrm{a}\nabla u, u, u|_{\mathcal{F}_h}),
\]
 where $\mathrm{a}:=\mathrm{c}^{-1}$, in the finite element space ${\bld{V}_h}\times {W_h}\times{M_h}$ where
\begin{alignat*}{3}
{\bld{V}}_h&:=\{\bld{v}\in \bld{L}^2(\Omega):\quad &&\bld{v}|_K\in
\bld{V}(K)\quad&&\forall K\in\Oh\},
\\
W_h&:=\{w\in {L}^2(\Omega):\quad &&w|_K\in
W(K)\quad&&\forall K\in\Oh\},
\\
M_h&:=\{\mu\in L^2(\mathcal{F}_h):\quad &&\mu|_F\in
M(F)\quad&&\forall F\in\mathcal{F}_h\}.
\end{alignat*}
The approximation is defined  as the solution of
\begin{alignat*}{1}
(\mathrm{c}\,{{\bld q}_h},{\bld v})_{\Oh} 
- ({{u_h}},\nabla \cdot {\bld v})_{\Oh} 
+ \langle {\uhat_h},{\bld v} \cdot {\bld n}\rangle_{\dOh} & = 0,\\
-({{\bld q}_h},\nabla w)_{\Oh} 
+ \langle {\qhat_h}\cdot {\bld n},w\rangle_{\dOh} &= (f,w)_{\Oh},\\
\langle\mu, {\qhat_h}\cdot\bld{n}\rangle_{\dOh\setminus\partial\Omega}&=0,\\
\langle\mu,{\uhat_h}\rangle_{\partial\Omega}&=\langle\mu, u_{D}\rangle_{\partial\Omega},
\end{alignat*}
for all $(\bld{v},w,\mu)\in{\bld{V}_h}\times {W_h}\times{M_h}$, and
\[
{\qhat_h}\cdot{\bld n}
 :={{\bld q}_h}\cdot\bld{n} +\tau( {{u_h}}-{\uhat_h})\bld \qquad\mbox{ on } \dOh.
\]
The function $\tau$ is called the stabilization function as it ensures the stability and well-posedness of the approximation. This completes the definition of the HDG methods. 

A few, brief comments are in order:
\begin{itemize}
\item  The  different HDG methods are
obtained by choosing the local spaces $\bld{V}(K)$,
$W(K)$, $M(F)$, and the stabilization function $\tau$.  The function $\tau$ is called the stabilization function as it ensures the stability and well-posedness of the approximation provided by the methods. 

\item  In \cite{CockburnGopalakrishnanLazarov09}, it was 
shown that, under suitable choices of these spaces and stabilization function, the resulting HDG 
methods are actually previously known DG  methods, none of which belong to the class of DG 
methods considered in the unifying analysis proposed in \cite{ArnoldBrezziCockburnMarini02}. 
The HDG methods are thus DG methods.

\item  The distinctive feature of the HDG methods is that their numerical traces are defined in such a way that the only globally-coupled degrees of freedom 
are  those of the numerical trace $\uhat_h$. This is why these methods are called {\em hybridizable}.

\item In \cite{CockburnGopalakrishnanLazarov09}, it was also 
pointed out that, when the stabilization function $\tau$ is identically zero, the {\em hybridized} form of the 
well known mixed methods is obtained.
\end{itemize}

\subsubsection{A second formulation of  the HDG methods}
To continue our discussion, we introduce two different, but strongly related, formulations of our model problem based in the following equivalent rewritings 
\begin{center}
\begin{minipage}[b]{0.46\textwidth}
\begin{equation*}
(F_{\mathrm{c}})\left\{ 
\begin{array}{rcll}
\boldsymbol{g}&=& \nabla u         &\mbox{ in } \Omega,\\
\mathrm{c}\,\boldsymbol{q}  &=& - \boldsymbol{g}               &\mbox{ in } \Omega,\\
\nabla\cdot \boldsymbol{q} &=& f              &\mbox{ in } \Omega,\\
u               &=& u_{D}          &\mbox{ on } \partial \Omega,  
\end{array}
\right.
\end{equation*}
\end{minipage}
\hfill  \hfill
\begin{minipage}[b]{0.46\textwidth}
\begin{equation*}
(F_{\mathrm{a}})\left\{ 
\begin{array}{rcll}
\boldsymbol{g}   &=& \nabla u           &\mbox{ in } \Omega,\\
\boldsymbol{q}   &=& -\mathrm{a}\,\boldsymbol{g}               &\mbox{ in } \Omega,\\
\nabla\cdot \boldsymbol{q} &=& f              &\mbox{ in } \Omega,\\
u               &=& u_{D}          &\mbox{ on } \partial \Omega,
\end{array}
\right.
\end{equation*}
\end{minipage}
\end{center}
where $\mathrm{a}:=\mathrm{c}^{-1}$. It has been noted long time ago that, if $\mathrm{c}$ is difficult to obtain from $\mathrm{a}$, which is usually the observed data, a numerical method based on the set of equations ($F_{\mathrm{a}}$) is in order. 

\subsubsection*{HDG methods based on ($F_{\mathrm{c}}$) and on ($F_{\mathrm{a}}$)}
Let us show that the HDG method displayed in the previous section is naturally associated with the set of equations ($F_{\mathrm{c}}$). Indeed, if we define the approximation
$({\bld{q}_h}, {\bld{g}_h}, {{u_h}},{\uhat_h})$,  
in the finite element space ${\bld{V}_h}\times{\bld{V}_h}\times {W_h}\times{M_h}$ 
as the solution of
\begin{alignat*}{1}
-({{\bld g}_h},{\bld v})_{\Oh} 
- ({{u_h}},\nabla \cdot {\bld v})_{\Oh} 
+ \langle {\uhat_h},{\bld v} \cdot {\bld n}\rangle_{\dOh} & = 0,\\
(\mathrm{c}\,{{\bld q}_h},{\bld v})_{\Oh} & = -({{\bld g}_h},{\bld v})_{\Oh} ,\\
-({{\bld q}_h},\nabla w)_{\Oh} 
+ \langle {\qhat_h}\cdot {\bld n},w\rangle_{\dOh} &= (f,w)_{\Oh},\\
\langle\mu, {\qhat_h}\cdot\bld{n}\rangle_{\dOh\setminus\partial\Omega}&=0,\\
\langle\mu,{\uhat_h}\rangle_{\partial\Omega}&=\langle\mu, u_{D}\rangle_{\partial\Omega},
\end{alignat*}
for all $(\bld{v},w,\mu)\in{\bld{V}_h}\times {W_h}\times{M_h}$, and
\[
{\qhat_h}\cdot{\bld n}
 :={{\bld q}_h}\cdot\bld{n} +\tau( {{u_h}}-{\uhat_h})\bld \qquad\mbox{ on } \dOh,
\]
we immediately see that the component $({\bld{q}_h}, {{u_h}},{\uhat_h})$ is nothing but 
the  approximate solution provided by the HDG method just considered, that is, the HDG method 
introduced in 2009 in \cite{CockburnGopalakrishnanLazarov09}.

Let us now introduce the HDG method based on the equations ($F_{\mathrm{a}}$). The HDG method defines the approximation
$({\bld{q}_h}, {\bld{g}_h}, {{u_h}},{\uhat_h})$,  
in the finite element space ${\bld{V}_h}\times{\bld{V}_h}\times {W_h}\times{M_h}$ 
as the solution of
\begin{alignat*}{1}
-({{\bld g}_h},{\bld v})_{\Oh} 
- ({{u_h}},\nabla \cdot {\bld v})_{\Oh} 
+ \langle {\uhat_h},{\bld v} \cdot {\bld n}\rangle_{\dOh} & = 0,\\
({{\bld q}_h},{\bld v})_{\Oh} & = - (\mathrm{a}\,{{\bld g}_h},{\bld v})_{\Oh},\\
-({{\bld q}_h},\nabla w)_{\Oh} 
+ \langle {\qhat_h}\cdot {\bld n},w\rangle_{\dOh} &= (f,w)_{\Oh},\\
\langle\mu, {\qhat_h}\cdot\bld{n}\rangle_{\dOh\setminus\partial\Omega}&=0,\\
\langle\mu,{\uhat_h}\rangle_{\partial\Omega}&=\langle\mu, u_{D}\rangle_{\partial\Omega},
\end{alignat*}
for all $(\bld{v},w,\mu)\in{\bld{V}_h}\times {W_h}\times{M_h}$, and
\[
{\qhat_h}\cdot{\bld n}
 :={{\bld q}_h}\cdot\bld{n} +\tau( {{u_h}}-{\uhat_h})\bld \qquad\mbox{ on } \dOh.
\]

Some remarks are in order:

\begin{itemize}
\item It is not difficult to see that for any HDG method defined by using the equations ($F_{\mathrm{c}}$), we can immediately define the corresponding HDG method for the equations 
($F_{\mathrm{a}}$). A detailed comparative study of HDG methods based on the equations ($F_{\mathrm{c}}$) and ($F_{\mathrm{a}}$) has been recently carried out in \cite{CockburnSanchezXiong18}.
There, conditions on the local spaces under which the approximations provided by both HDG methods are superclose are provided.

\item HDG methods based on the equations of the type ($F_{\mathrm{a}}$) were introduced 
for the more difficult problem of linear (and nonlinear) elasticity. In those equations,
{\color{black}$\bld{g}$ corresponds to the strain, $\bld{q}$ to the stress, $\mathrm{a}$ to the so-called {\em constitutive} tensor and $\mathrm{c}\,$ to} the so-called {\em compliance} tensor. 
The first HDG method based on the equations corresponding to ($F_{\mathrm{a}}$)
for linear and nonlinear elasticity were introduced in the 2008 PhD. thesis \cite{SoonThesis08}
part of which was then published in 2009 in \cite{SoonCockburnStolarski09}. In 2015, the 
convergence properties of the HDG methods were re-examined in \cite{FuCockburnStolarski15} for 
the linear case and  in \cite{KabariaLewCockburn15} for the nonlinear case. An HDG method  based on the equations ($F_{\mathrm{a}}$) was introduced for the $p$-Laplacian in 
2016 in \cite{CockburnShen16}. 

\item The HDG methods based on the equations ($F_{\mathrm{a}}$) are also {\em hybridizable}.

\item As with the HDG methods based on the equations ($F_{\mathrm{c}}$), when the stabilization function $\tau$ can be taken to be equal to zero, we obtain mixed methods. These mixed methods
were explored in the mid 90's for the first time in \cite{Koebbe93,Chen94,ArbogastWheelerYotov97}. 
If we formulate the 1985 mixed BDM method
\cite{BrezziDouglasMarini85}, defined by using the equations ($F_{\mathrm{c}}$), by using the equations ($F_{\mathrm{a}}$), we obtain the 1994 mixed method proposed in \cite{Chen94} -- for a nonlinear diffusion equation. The author found that these two methods were so similar, he still refered to the resulting method as a BDM mixed method. He did not find the need to change the name of the method. An extension of this work to arbitrary mixed methods was carried out in 1997  \cite{ArbogastWheelerYotov97}. Therein, the authors called the resulting mixed methods the ``expanded'' mixed methods, in order to stress the fact that a new unknown for the gradient had to be added to the equations. When these methods are hybridized, and their space of gradients $\widetilde{\boldsymbol{V}}_h$ is such that $\widetilde{\boldsymbol{V}}_h\cap \boldsymbol{H}(div,\Omega)= \boldsymbol{V}_h$, their formulation fits the HDG formulation based on ($F_{\mathrm{a}}$) with a zero stabilization function.

\end{itemize}

\subsection{The 2013 WG methods \cite{WangYeMain13} are mixed methods}
\label{subsection:firstWGmethods}
Let us recall the fact that, the very first WG methods \cite{WangYeMain13} were devised to discretize the steady-state
second-order elliptic equation
\[
-\nabla\cdot (a\nabla u)+ \nabla\cdot(bu)+ c u=f.
\]
Since the main contribution in  \cite{WangYeMain13} is the discretization of the second-order term, the main result we present here is for the case in which  $b=0$ and $c=0$.  

So, for $b=0$ and $c=0$, we are then going to show that the very first WG methods 
\cite{WangYeMain13} are nothing but a rewriting of the hybridized version of the mixed methods introduced in 1994 \cite{Chen94} and in 1997 \cite{ArbogastWheelerYotov97}. 

Since these two methods satisfy the formulation of the HDG method based on the equations ($F_{\mathrm{a}}$) with a zero stabilization function $\tau$, to prove this, we only have to rewrite such formulation 
 in terms of $(u_h, \uhat_h)$.  To achieve this, we express
the approximate gradient $\bld{g}_h$ and the approximate flux $\bld{q}_h$ as linear mappings of 
$\mathsf{u}_h:=({u_h},{\uhat_h})$ defined by using the first and second equations, respectively, 
defining the HDG method. 

Thus, for any $\mathsf{w}_h:=(w,\what)\in L^2(\Oh)\times L^2(\dOh)$, we define $({{\bld{G}_{\mathsf{w}_h}}},{{\bld{Q}_{\mathsf{w}_h}}})\in \bld{V}_h\times \bld{V}_h$ to be the solution of
\begin{alignat*}{2}
-({{\bld{G}_{\mathsf{w}_h}}},{\bld v})_{\Oh} 
&= (w,\nabla \cdot {\bld v})_{\Oh} 
 - \langle {\what},{\bld v} \cdot {\bld n}\rangle_{\dOh}
&&\quad\forall  {\bld{v}}\in \bld{V}_h,
\\
({{\bld{Q}_{\mathsf{w}_h}}},{\bld v})_{\Oh} & = - (\mathrm{a}\,{{\bld{G}_{\mathsf{w}_h}}},{\bld v})_{\Oh}
&&\quad\forall  {\bld{v}}\in \bld{V}_h.
\end{alignat*}
Note that, when $\mathsf{w}_h:=\mathsf{u}_h=(u_h,\uhat_h)$, we have that 
$\bld{g}_h={{\bld{G}_{\mathsf{u}_h}}}$
by the first equation defining the HDG methods, and $\bld{q}_h={{\bld{Q}_{\mathsf{u}_h}}}$ by the 
second. Next we show how, by using these mappings, we can eliminate $\bld{g}_h$ and $\bld{q}_h$ from the 
equations. 

Integrating by parts the third equation defining the HDG methods, we get
\begin{alignat*}{1}
(\nabla\cdot {{\bld q}_h},w)_{\Oh}
+\langle (\qhat_h-{\bld q}_h)\cdot {\bld n},w\rangle_{\dOh} = (f,w)_{\Oh}\quad\forall w\in W_h.
\end{alignat*} 
Taking $\bld{v}:={\bld q}_h$ in the equation defining ${{\bld{G}_{\mathsf{w}_h}}}$, we get that
\begin{alignat*}{1}
-({{\bld{G}_{\mathsf{w}_h}}},{\bld {q}}_h)_{\Oh} 
= (w,\nabla \cdot {\bld q}_h)_{\Oh} 
 - \langle {\what},{\bld q}_h \cdot {\bld n}\rangle_{\dOh}.
\end{alignat*}
If in the equation defining ${{\bld{Q}_{\mathsf{w}_h}}}$, we set 
$\mathsf{w}_h:=\mathsf{u}_h$ and the take
$\bld{v}:={{\bld{G}_{\mathsf{w}_h}}}$, we get that
\[
(\bld{q}_h, {{\bld{G}_{\mathsf{w}_h}}})=({{\bld{Q}_{\mathsf{u}_h}}},{{\bld{G}_{\mathsf{w}_h}}})_{\Oh} = - (\mathrm{a}\,{{\bld{G}_{\mathsf{u}_h}}},{{\bld{G}_{\mathsf{w}_h}}})_{\Oh}.
\]
This implies that
\begin{alignat*}{1}
(\mathrm{a}\,{{\bld{G}_{\mathsf{w}_h}}}, {{\bld{G}_{\mathsf{u}_h}}})_{\Oh} +\langle {\what},{\bld q}_h \cdot {\bld n}\rangle_{\dOh}
+\langle (\qhat_h-{\bld q}_h)\cdot {\bld n},w\rangle_{\dOh} = (f,w)_{\Oh}\quad\forall w\in W_h.
\end{alignat*} 
Adding and subtracting the term $\langle {w},\qhat_h \cdot {\bld n}\rangle_{\dOh}$, and rearranging \ terms, we get
\begin{alignat*}{1}
(\mathrm{a}\,{{\bld{G}_{\mathsf{w}_h}}}, {{\bld{G}_{\mathsf{u}_h}}})_{\Oh}
+\langle (\qhat_h-{\bld q}_h)\cdot {\bld n},w-\what \rangle_{\dOh} &= (f,w)_{\Oh}
-\langle {\what},\qhat_h \cdot {\bld n}\rangle_{\dOh}
\\
&= (f,w)_{\Oh}
-\langle {\what},\qhat_h \cdot {\bld n}\rangle_{\partial\Omega},
\end{alignat*} 
for all $(w,\what)\in W_h\times M_h$,
by the fourth equation defining the HDG methods. Now, inserting the definition of the numerical trace
$\qhat_h \cdot {\bld n}$ and noting that $\bld{q}_h={{\bld{Q}_{\mathsf{u}_h}}}$, we get
\begin{alignat*}{1}
(\mathrm{a}\,{{\bld{G}_{\mathsf{w}_h}}}, {{\bld{G}_{\mathsf{u}_h}}})_{\Oh}
+\langle \tau(u_h-\uhat_h), w-\what \rangle_{\dOh}&= (f,w)_{\Oh}
-\langle {\what},\qhat_h \cdot {\bld n}\rangle_{\partial\Omega},
\end{alignat*} 
for all $(w,\what)\in W_h\times M_h$.

Finally, if we set
$
M_h(\eta):=\{\zeta\in M_h: \langle\zeta,\mu\rangle_{\partial\Omega}= \langle\eta,\mu\rangle_{\partial\Omega}\quad\forall \mu\in M_h\},
$
the formulation of the HDG method determines the approximation to the exact solution
$(u|_\Omega,u|_{\mathcal{F}_h})$, $\mathsf{u}_h:=(u_h,\uhat_h)$, as the element of $W_h\times M_h(g)$ such that
\begin{alignat*}{1}
(\mathrm{a}\,{{\bld{G}_{\mathsf{w}_h}}}, {{\bld{G}_{\mathsf{u}_h}}})_{\Oh}
+\langle \tau(u_h-\uhat_h), w-\what \rangle_{\dOh}&= (f,w)_{\Oh},
\end{alignat*}
for all $\mathsf{w}_h:=(w,\what)\in W_h\times M_h$.

These are precisely the WG methods in \cite{WangYeMain13} obtained by taking
the stabilization function $\tau=0$, the local space $M(F):=\mathcal{P}_k(F)$, 
and
\begin{alignat*}{3}
&\bld{V}(K) := [\mathcal{P}_k(K)]^d,
&&\quad
W(K):=\mathcal{P}_{k-1}(K),
&&\quad 
\text{(Example 5.1)},
\\
&\bld{V}(K) := [\mathcal{P}_k(K)]^d + {\bld x} \mathcal{P}_k(K),
&&\quad
W(K):=\mathcal{P}_{k}(K),
&&\quad
\text{(Example 5.2)}.
\end{alignat*}
Since, the stabilization function $\tau$ is zero, these are mixed methods. In fact the first is nothing but 
the hybridized version of the BDM method considered in 1994 \cite{Chen94}, and the second is 
nothing but the hybridized version of the RT method considered in 1997 
\cite{ArbogastWheelerYotov97}. Thus, for $b=0$ and $c=0$, the first WG methods \cite{WangYeMain13} are nothing but the 
hybridized version of mixed methods.

\subsubsection*{Remark}
Let us emphasize that taking $c\ne0$ does not alter the above results.
Concerning the discretization of the convective term (when $b\ne0$), the way it was done in \cite{WangYeMain13} represents a regression with respect to the original upwinding technique introduced back in the early 70's in \cite{ReedHill73,LesaintRaviart74} since the absence of upwinding introduces a lack of stabilization.  This technique was {\bf rescued} a few years later in \cite{ChenFengXie17}.

The original definition of the WG methods, equation (4.6) in \cite{WangYeMain13} holds for piecewise polynomial approximations (defined on meshes made of triangles only), but the methods were proven to be well posed only for the spaces of the
Examples 5.1 and 5.2, which constitute classic mixed methods, as just proved above. If the spaces defining the  WG methods were not associated to classic mixed methods, the brief attempt of relating the WG and the HDG methods
 made by the authors, see our Introduction, would then be granted. However, even in this case, the stabilization function the authors use is zero and so, the methods are not HDG methods but only mixed methods, as pointed out back in 2009 \cite{CockburnGopalakrishnanLazarov09}. {\em However}, if $c=0$, then the reaction term does not play a stabilizing role anymore. Hence, a nonzero 
 stabilization function $\tau$ would have to be introduced, and the method would become an HDG method. This is exactly what happened with the WG methods we consider next.

\subsection{Equivalence of the HDG and WG methods}

\subsubsection{The 2014 WG method \cite{WangYemixedsecondorder14} is an HDG method}
Here, we establish the equivalence between the 2014 WG method proposed 
in \cite{WangYemixedsecondorder14} and the HDG method. To do that, 
we need to pick a particular stabilization function $\tau$, and then express the HDG method (using the first formulation) as a method whose 
``unknowns"  are $({\bld{q}_h}, {{u_h}},{\qhat_h})$ are instead of $({\bld{q}_h}, {{u_h}},{\uhat_h})$. 

To do this, we rely on two very simple observations. The first is that, if we take as stabilization function $\tau$, the so-called Lehrenfeld-Sch\"oberl stabilization \cite{Lehrenfeld10,Schoeberl09} function, that is,
\[
\tau:=\frac{\rho}{h_K} P_{M(F)} \quad\text{ on the face $F$ of the element $K\in \Oh$,} 
\]
where $h_K$ is the diameter of the element $K$ and $P_{M(F)}$ is the
$L^2$-projection into $M(F)$, we have that
\[
\qhat_h\cdot\bld{n}= \bld{q}_h\cdot\bld{n}+ \frac{1}{h} (P_M(u_h) -\uhat_h),
\]
which immediately implies that
\[
\uhat_h= P_M(u_h)+ \frac{h}{\rho} (\bld{q}_h-\qhat_h)\cdot\bld{n}.
\]
The second observation is that, for each face $F$ of the mesh, since $\uhat_h|_F\in M(F)$, 
when
 \[
 \bld{V}(K)\cdot\bld{n}|_F\subset M(F),
 \]
 we have that $\qhat_h\cdot\bld{n}|_F$ also belongs to 
$M(F)$. In other words, by the third equation defining the HDG method, the numerical trace $\qhat{}_h$ belongs to the space
\[
\bld{N}_h:=\{\boldsymbol{\nu}\in\bld{L}^2(\mathcal{F}_h):\boldsymbol{\nu}\cdot\bld{n}|_{F}\in M(F)\forall F\in \mathcal{F}_h\},
\]
where the superscripts $\pm$ denote the corresponding traces.

We can thus rewrite the HDG method as follows.
The approximate solution given by the HDG method is the function 
$(\bld{q}_h,u_h,\qhat_h)\in \bld{V}_h\times W_h\times \bld{N}_h$ satisfying the equations
\begin{alignat*}{2}
(\mathrm{c}\,\bld{q}_h,{\bld v})_{\Oh} 
- (u_h,\nabla \cdot {\bld v})_{\Oh} 
+ \langle \uhat_h,{\bld v} \cdot {\bld n}\rangle_{\dOh} & = 0
&&\quad\forall \bld{v}\in\bld{V}_h,\\
-(\bld{q}_h,\nabla w)_{\Oh} 
+ \langle \qhat_h \cdot {\bld n},w\rangle_{\dOh} &= ({f},w)_{\Oh}
&&\quad\forall w\in W_h,
\\
\uhat_h= P_M(u_h)+ \frac{h}{\rho} (\bld{q}_h-\qhat_h)\cdot\bld{n},
&&&\quad\mbox{on }\dOh,
\\
\langle\uhat_h,\boldsymbol{\nu}\cdot\bld{n}\rangle_{\dOh} &=\langle u_D,\boldsymbol{\nu}\cdot\bld{n}\rangle_{\partial \Omega}
&&\quad\forall \boldsymbol{\nu}\in \bld{N}_h.
\end{alignat*}

Now, let us work on the first equation. If we subtract the fourth equation from the first one, we immediately get that
{\small
\begin{alignat*}{1}
-\langle u_D,\boldsymbol{\nu}\cdot\bld{n}\rangle_{\partial \Omega}
=
&
(\mathrm{c}\,\bld{q}_h,{\bld v})_{\Oh} 
- (u_h,\nabla \cdot {\bld v})_{\Oh} 
+ \langle \uhat_h,({\bld v}-{\bld\nu})  \cdot {\bld n}\rangle_{\dOh}
\\
=
&(\mathrm{c}\,\bld{q}_h,{\bld v})_{\Oh} 
\!+\! (\nabla u_h, {\bld v})_{\Oh} - \langle u_h, {\bld v}  \cdot {\bld n}\rangle_{\dOh}
\!+\! \langle \uhat_h,({\bld v}-{\bld\nu})  \cdot {\bld n}\rangle_{\dOh}
\\
=
&
(\mathrm{c}\,\bld{q}_h,{\bld v})_{\Oh} 
\!+\! (\nabla u_h, {\bld v})_{\Oh} - \langle u_h, {\bld 
\nu}  \cdot {\bld n}\rangle_{\dOh}
\!+\! \langle \uhat_h-u_h,({\bld v}-{\bld\nu})  \cdot {\bld n}\rangle_{\dOh}
\\
=
&
(\mathrm{c}\,\bld{q}_h,{\bld v})_{\Oh} 
\!+\! (\nabla u_h, {\bld v})_{\Oh} - \langle u_h, {\bld 
\nu}  \cdot {\bld n}\rangle_{\dOh}
\\
&+\! \langle \frac{h}{\rho} (\bld{q}_h-\qhat_h)\cdot\bld{n},({\bld v}-{\bld\nu})  \cdot {\bld n}\rangle_{\dOh},
\end{alignat*}
}
since  $ \bld{V}(K)\cdot\bld{n}|_F\subset M(F)$ and since $\uhat_h= P_M(u_h)+ \frac{h}{\rho} (\bld{q}_h-\qhat_h)\cdot\bld{n}$ on $\dOh$.

Thus, the approximate solution given by the HDG method is the function 
 satisfying the equations
\begin{alignat*}{4}
&[(\mathrm{c}\,\bld{q}_h,{\bld v})_{\Oh} + \langle \frac{h}{\rho} (\bld{q}_h-\qhat_h)\cdot\bld{n}, &&({\bld v}-\bld{\nu}) \cdot {\bld n}\rangle_{\dOh}]
\\
&&&- [-(\nabla u_h,{\bld v})_{\Oh}  + \langle u_h, {\bld 
\nu}  \cdot {\bld n}\rangle_{\dOh} ] && = 0,\\
&&&\;\,\phantom{-}[-(\bld{q}_h,\nabla w)_{\Oh} 
+ \langle \qhat_h \cdot {\bld n},w\rangle_{\dOh} ]&&= ({f},w)_{\Oh},
\end{alignat*}
for all $(\bld{v},\omega,\boldsymbol{\nu})\in \bld{V}_h\times W_h\times \bld{N}_h$.

For the choice of local spaces $\bld{V}(K) := [\mathcal{P}_k(K)]^d$,
$W(K):=\mathcal{P}_{k+1}(K)$,
$M(F):=\mathcal{P}_k(F)$,
this is precisely the 2014 WG method proposed in \cite{WangYemixedsecondorder14}; see equations (3.14) and (3.15) therein. In other words, the 2014 WG method in \cite{WangYemixedsecondorder14} is nothing but an HDG method.

\subsubsection{The 2015  WG method in \cite{MuWangYepolytopal15} is an HDG method}

To show that the 2015  WG method in \cite{MuWangYepolytopal15} is an HDG method, we simply rewrite the formulation of the HDG method based on ($F_{\mathrm{a}}$) in terms of $(u_h, \uhat_h)$ only, as we did in Subsection \ref{subsection:firstWGmethods}.

Thus, with the same notation used there, we obtain that the HDG method determines the approximation to the exact solution
$(u|_\Omega,u|_{\mathcal{F}_h})$, $\mathsf{u}_h:=(u_h,\uhat_h)$, as the element of $W_h\times M_h(u_D)$ such that
\begin{alignat*}{1}
(\mathrm{a}\,{{\bld{G}_{\mathsf{w}_h}}}, {{\bld{G}_{\mathsf{u}_h}}})_{\Oh}
+\langle \tau(u_h-\uhat_h), w-\what \rangle_{\dOh}&= (f,w)_{\Oh},
\end{alignat*}
for all $\mathsf{w}_h:=(w,\what)\in W_h\times M_h$. 

For the choice of spaces $\bld{V}(K) := [\mathcal{P}_k(K)]^d$,
$W(K):=\mathcal{P}_{k+1}(K)$,
$M(F):=\mathcal{P}_{k+1}(F)$,
and the stabilization function 
\[
\tau:= \frac{\rho}{h_K} \quad\text{ on the face $F$ of the element $K\in \Oh$,}
\]
where $\rho$ is a positive constant, 
this is precisely the WG method proposed in 2015 \cite{MuWangYepolytopal15}. Thus,  the 2015 WG method in \cite{MuWangYepolytopal15} is nothing but an HDG method.

\subsubsection{The 2015 WG method in \cite{MuWangYeLS15} is an HDG method}
The 2015 WG method in \cite{MuWangYeLS15} is obtained by taking
$\bld{V}(K) := [\mathcal{P}_k(K)]^d$,
$W(K):=\mathcal{P}_{k+1}(K)$,
$M(F):=\mathcal{P}_k(F)$,
and the stabilization function is 
\[
\tau:= \frac{\rho}{h_K}P_{M(F)} \quad\text{ on the face $F$ of the element $K\in \Oh$,}
\]
where $\rho$ is a positive constant.

Thus, the only differences between the 2015 WG method proposed in \cite{MuWangYeLS15}  and the WG method proposed in 2015 \cite{MuWangYepolytopal15} are the choice of the local space $M(F)$ and of the related
stabilization function $\tau$. For the WG method proposed in 2015 \cite{MuWangYepolytopal15}, $M(F):=\mathcal{P}_{k+1}(F)$ and the stabilization function is a simple multiplication function, 
whereas  for the 2015 WG method in \cite{MuWangYeLS15}, $M(F):=\mathcal{P}_{k}(F)$ and the stabilization function is the Lehrenfeld-Sch\"oberl stabilization function \cite{Lehrenfeld10,Schoeberl09}. Since we proved in the previous subsection that the WG method proposed in 2015 \cite{MuWangYepolytopal15} is an HDG method, the same holds for the  2015 WG method in \cite{MuWangYeLS15}. Thus, the 2015 WG method in \cite{MuWangYeLS15} is nothing but an HDG method.

Let us end by pointing out that the 2014 WG method in \cite{WangYemixedsecondorder14} is {\em identical} to the WG method in \cite{MuWangYeLS15} when the tensors $\mathrm{c}\,$ and $\mathrm{a}$ are piecewise constant. Even if these functions are not piecewise  constant, 
by the results of \cite{CockburnSanchezXiong18}, these two methods are {\em superclose}.

\section{The Stokes equations} 
The WG method for Stokes obtained in \cite{WangYeStokes16}
can be shown to be an HDG method with the local spaces 
\[
\mathrm{G}(K) := [\mathcal{P}_k(K)]^{d\times d},\;
{\bld V}(K):=[\mathcal{P}_{k+1}(K)]^d,\;
P(K):=\mathcal{P}_{k}(K),\;
{\bld M}(F):=[\mathcal{P}_k(F)]^d,
\]
and the Leherenfeld-Sch\"oberl stabilization function \cite{Lehrenfeld10,Schoeberl09}
\[
\tau:= \frac{1}{h_K}P_{M(F)} \quad\text{ on the face $F$ of the element $K\in \Oh$,}
\]
where $h_K$ is the diameter of the element $K$ and $P_{M(F)}$ is the
$L^2$-projection into the space ${\bld M}(F)$. 

\section{The biharmonic}
Here we show that the WG methods in 
\cite{MuWangWangYebiharmonic13,MuWangYebiharmonic14,WangWangbiharmonic14} 
for the biharmonic problem
\begin{alignat*}{2}
\Delta^2 u  = f \quad
\text{ on } \Omega,
\quad
u  = \frac{\partial u}{\partial\boldsymbol{n}}=0\quad \text{ on } \partial\Omega,
\end{alignat*}
are also rewritings of the HDG methods. To do that, it is enough to follow 
a straightforward extension of the procedure used in the previous cases.

\subsection{The HDG methods}

To define an HDG method,  we first 
rewrite the problem as a first-order system. Here, we consider two ways of doing that.
\subsubsection{A first formulation for HDG methods}
The first way uses the following equations:
\begin{alignat*}{2}
{\boldsymbol{\sigma}}= \nabla\cdot z,
\quad
\nabla \cdot{\boldsymbol{\sigma}} = f,
\quad
{\boldsymbol{q}} =\nabla u,
\quad
\nabla\cdot {\boldsymbol{q}}  = z
\quad    \text { in } \Omega,
\end{alignat*}
with the boundary conditions
$
  u = 0
  \text{ and }
  {\boldsymbol{q}}\cdot{\boldsymbol{n}} =0  \text{ on } \dO.
$

The HDG method seeks an approximation $({\boldsymbol{\sigma}}_h,z_h,{\boldsymbol{q}}_h,u_h)$ to the exact solution
\[
(\boldsymbol{\sigma}, z, \boldsymbol{q},u) \quad\text{ in } \Omega,
\]
in a finite dimensional space 
${\boldsymbol{\Sigma}}_h\times Z_h\times {\boldsymbol{Q}}_h\times W_h$ of the form
\begin{alignat*}{2}
{\boldsymbol{\Sigma}}_h:=&\lbrace {\boldsymbol{v}}\in {\boldsymbol{L}^2}(\Omega) \,: &&\quad\, {\boldsymbol{v}}|_K \in
     {\boldsymbol{\Sigma}}(K)\quad\forall K\in\Oh\rbrace,\\
{Z}_h:=&\lbrace s\in {L}^2(\Omega) \,: &&\quad\, s|_K \in
     Z(K)\quad\forall K\in\Oh\rbrace,\\
{\boldsymbol{Q}}_h:=&\lbrace {\boldsymbol{m}}\in {\boldsymbol{L}^2}(\Omega) \,: 
&&\quad
\, {\boldsymbol{m}}|_K \in
     {\boldsymbol{Q}}(K)\quad\forall K\in\Oh\rbrace,\\
W_h:=&\lbrace \omega \in {L^2}(\Omega)\,:
&&\quad\, \omega|_K \in
     W(K)\quad\forall K\in\Oh\rbrace,
\end{alignat*}
and an approximation $(\widehat{\boldsymbol{\sigma}}_h\cdot\boldsymbol{n}, \widehat{z}_h, \widehat{\boldsymbol{q}}_h\cdot\boldsymbol{n},, \widehat{u}_h)$ to the traces of the exact solution
\[
(\boldsymbol{\sigma}\cdot\boldsymbol{n}, z, \boldsymbol{q}\cdot\boldsymbol{n},u) \quad\text{ on } \partial\mathcal{T}_h,
\]
and determines them by requiring that
\begin{alignat*}{1}
({\boldsymbol{\sigma}}_h, {\boldsymbol{m}})_\Oh+(z_h,\nabla{\boldsymbol{m}})_\Oh-\langle
\widehat{z}_h,{\boldsymbol{m}}\cdot {\boldsymbol{n}}\rangle_\dOh=&0,
\\
-({\boldsymbol{\sigma}}_h,\nabla\omega)_\Oh+\langle \widehat{{\boldsymbol{\sigma}}}_h\cdot{\boldsymbol{n}},
\omega\rangle_\dOh=&(f,\omega)_\Oh,
\\
({\boldsymbol{q}}_h, {\boldsymbol{v}})_{\Oh}+(u_h,\nabla\cdot{\boldsymbol{v}})_{\Oh}-\langle
\widehat{u}_h,{\boldsymbol{v}}\cdot{\boldsymbol{n}} \rangle_{\dOh}=&0,
\\
-({\boldsymbol{q}}_h, \nabla s)_{\Oh}+\langle
\widehat{{\boldsymbol{q}}}_h\cdot{\boldsymbol{n},s}\rangle_{\dOh}=&(z_h,s)_{\Oh},
\end{alignat*}
for all $({\boldsymbol{v}},s,{\boldsymbol{m}},\omega)\in {\boldsymbol{\Sigma}}_h\times Z_h\times {\boldsymbol{Q}}_h\times W_h$. 
As is typical of all the HDG methods, to complete their definition, we must impose the boundary conditions and define the approximate traces in such a way that the approximation is well defined and that the method can be {\em statically condensed}. 

A couple of remarks are in order:

\begin{itemize}

\item All the HDG methods based on this formulation are 
obtained by varying the local spaces and the definition of the approximate traces. 

\item The very first HDG method of the biharmonic, introduced in 2009 \cite{CockburnDongGuzmanBiharmonic09}, is obtained when we take the local spaces
to be
\[
\boldsymbol{\Sigma}(K)=\boldsymbol{Q}(K):=[{\mathcal{P}}_k(K)]^d
\quad
\text{ and }\quad
Z(K)=W(K):=\mathcal{P}_k(K),
\]
impose the boundary conditions
\[
\widehat{u}_h=0
\quad
\text{ and } 
\widehat{{\boldsymbol{q}}}_h\cdot{\boldsymbol{n}}=0
\qquad
\text{ on }\dO,
\]
pick the numerical traces $\widehat{z}_h$ and $\widehat{u}_h$  as the globally coupled {\em unknowns} in the finite dimensional space 
\[
M_h:=\{\mu\in L^2(\mathcal{F}_h): \mu|_F\in \mathcal{P}_k(F)\quad\forall \; F\in \mathcal{F}_h\},
\]
define the numerical traces $( \widehat{{\boldsymbol{q}}}_h\cdot\boldsymbol{n},  \widehat{{\boldsymbol{\sigma}}}_h\cdot\boldsymbol{n})$ by
\begin{alignat*}{2}
\widehat{{\boldsymbol{q}}}_h\cdot\boldsymbol{n} 
&=
{{\boldsymbol{q}}}_h\cdot\boldsymbol{n} -
\tau\;(u_h-\widehat{u}_h)
&&\qquad\mbox{ on } \partial \Oh\setminus\dO,
\\
\widehat{{\boldsymbol{\sigma}}}_h\cdot\boldsymbol{n}
&={{\boldsymbol{\sigma}}}_h\cdot\boldsymbol{n}-\tau\;(z_h-\widehat{z}_h),
&&\qquad\mbox{ on } \partial \Oh,
\end{alignat*}
where $\tau$ is the stabilization function of the SFH method for scalar diffusion introduced in 2008 \cite{CockburnDongGuzmanSFH08}, and impose the transmission conditions
\begin{alignat*}{2}
\langle{{\boldsymbol{q}}}_h\cdot\boldsymbol{n}, \widehat{\omega}\rangle_{\dOh\setminus\dO}&=0 &&\quad\forall \widehat{\omega}\in M_h,
\\
\langle{{\boldsymbol{\sigma}}}_h\cdot\boldsymbol{n}, \widehat{\omega}\rangle_{\dOh\setminus\dO}&=0 &&\quad\forall \widehat{\omega}\in M_h.
\end{alignat*}

\item The above HDG method extends the HDG methods for Timoshenko beams introduced in 2010 \cite{CelikerCockburnShiTimo10} and then analyzed in 2012 \cite{CelikerCockburnShiTimo12} . 

\end{itemize}

\subsubsection{A second formulation for HDG methods}

A second way of writing the biharmonic equations is
\begin{alignat*}{2}
{\boldsymbol{\sigma}}= \nabla\cdot \underline{\boldsymbol{z}},
\quad
\nabla \cdot{\boldsymbol{\sigma}}  = f,
\quad
{\boldsymbol{q}} =\nabla u,
\quad
\nabla {\boldsymbol{q}}  = \underline{\boldsymbol{z}}
\quad  \text { in } \Omega,
\end{alignat*}
with the boundary conditions
$
  u = 0
  \text{ and }
  {\boldsymbol{q}}\cdot{\boldsymbol{n}} =0  \text{ on } \dO.
$
Note that the only difference between these equations is that whereas $z$ is the Laplacian of $u$,  
$\underline{\boldsymbol{z}}$ is its Hessian.

The HDG method seeks an approximation $({\boldsymbol{\sigma}}_h,\underline{\boldsymbol{z}}_h,{\boldsymbol{q}}_h,u_h)$ to the exact solution
\[
(\boldsymbol{\sigma}, \underline{\boldsymbol{z}}, \boldsymbol{q},u) \quad\text{ in } \Omega,
\]
in a finite dimensional space 
${\boldsymbol{\Sigma}}_h\times \underline{\boldsymbol{Z}}_h\times {\boldsymbol{Q}}_h\times W_h$ where
\begin{alignat*}{2}
{\boldsymbol{\Sigma}}_h:=&\lbrace {\boldsymbol{v}}\in {\boldsymbol{L}^2}(\Omega) \,: &&\quad\, {\boldsymbol{v}}|_K \in
     {\boldsymbol{\Sigma}}(K)\quad\forall K\in\Oh\rbrace,\\
{\underline{\boldsymbol{Z}}}_h:=&\lbrace \underline{\boldsymbol{s}}\in [{L}^2(\Omega)]^{d\times d} \,: &&\quad\, \underline{{\boldsymbol{s}}}|_K \in
     \underline{\boldsymbol{Z}}(K)\quad\forall K\in\Oh\rbrace,\\
{\boldsymbol{Q}}_h:=&\lbrace {\boldsymbol{m}}\in {\boldsymbol{L}^2}(\Omega) \,: 
&&\quad
\, {\boldsymbol{m}}|_K \in
     {\boldsymbol{Q}}(K)\quad\forall K\in\Oh\rbrace,\\
W_h:=&\lbrace \omega \in {L^2}(\Omega)\,:
&&\quad\, \omega|_K \in
     W(K)\quad\forall K\in\Oh\rbrace,
\end{alignat*}
and an approximation $(\widehat{\boldsymbol{\sigma}}_h\cdot\boldsymbol{n}, \widehat{\underline{\boldsymbol{z}}}_h\boldsymbol{n}, \widehat{\boldsymbol{q}}_h, \widehat{u}_h)$ to the traces of the exact solution
\[
(\boldsymbol{\sigma}\cdot\boldsymbol{n}, \underline{\boldsymbol{z}}\boldsymbol{n}, \boldsymbol{q},u) \quad\text{ on } \partial\mathcal{T}_h,
\]
and determines them by requiring that
\begin{alignat*}{1}
({\boldsymbol{\sigma}}_h, {\boldsymbol{m}})_\Oh+(\underline{\boldsymbol{z}}_h,\nabla{\boldsymbol{m}})_\Oh-\langle
\widehat{\underline{\boldsymbol{z}}}_h {\boldsymbol{n}},{\boldsymbol{m}}\rangle_\dOh=&0,
\\
-({\boldsymbol{\sigma}}_h,\nabla\omega)_\Oh+\langle \widehat{{\boldsymbol{\sigma}}}_h\cdot{\boldsymbol{n}},
\omega\rangle_\dOh=&(f,\omega)_\Oh,
\\
({\boldsymbol{q}}_h, {\boldsymbol{v}})_{\Oh}+(u_h,\nabla\cdot{\boldsymbol{v}})_{\Oh}-\langle
\widehat{u}_h,{\boldsymbol{v}}\cdot{\boldsymbol{n}} \rangle_{\dOh}=&0,
\\
-({\boldsymbol{q}}_h, \nabla\cdot\underline{\boldsymbol{s}})_{\Oh}+\langle
\widehat{{\boldsymbol{q}}}_h,\underline{\boldsymbol{s}}{\boldsymbol{n}}\rangle_{\dOh}=&(\underline{\boldsymbol{z}}_h,\underline{\boldsymbol{s}})_{\Oh},
\end{alignat*}
for all $({\boldsymbol{v}},\underline{\boldsymbol{s}},{\boldsymbol{m}},\omega)\in {\boldsymbol{\Sigma}}_h\times \underline{\boldsymbol{Z}}_h\times {\boldsymbol{Q}}_h\times W_h$. Again, to complete the definition of the HDG methods, we must impose the boundary conditions and define the approximate traces in such a way that the approximation is well defined and that the method can be {\em statically condensed}. 

A couple of remarks are in order:

\begin{itemize}

\item All the HDG methods based on this formulation are 
obtained by varying the local spaces and the definition of the approximate traces.

\item A hybridized version of the mixed method for this formulation  introduced in 2011 \cite{BehrensGuzmanBiharmonic11} is obtained when we take the local spaces to be
\begin{alignat*}{1}
\boldsymbol{\Sigma}(K)&:=[\mathcal{P}_k(K)]^d\oplus \boldsymbol{x}\, \mathcal{P}_k(K),
\\
\underline{\boldsymbol{Z}}(K)&:=\{\underline{\boldsymbol{s}}\in [\mathcal{P}_k(K)]^{d\times d}: \;\text{each row of $\underline{\boldsymbol{s}}$ belongs to $\boldsymbol{\Sigma}(K)$}\},
\\
\boldsymbol{Q}(K)&:=[{\mathcal{P}}_k(K)]^d,
\\
W(K)&:=\mathcal{P}_k(K),
\end{alignat*}
impose the boundary conditions
\[
\widehat{u}_h =0,
\quad
\widehat{{\boldsymbol{q}}}_h=0\qquad
\text{ on } \dO,
\]
pick the traces $(\widehat{\boldsymbol{q}}_h, \widehat{u}_h)$ in the space $\boldsymbol{M}_h\times M_h$, where
\begin{alignat*}{3}
\boldsymbol{M}_h:=&\{\boldsymbol{\mu}\in [L^2(\mathcal{F}_h)]^d:&&\; \boldsymbol{\mu}|_K\in [\mathcal{P}_k(K)]^d&&\;\;\forall K\in \Oh\},
\\
{M}_h:=&\{{\mu}\in L^2(\mathcal{F}_h):&&\; {\mu}|_K\in \mathcal{P}_k(K)&&\;\;\forall K\in \Oh\},
\end{alignat*}
define the traces $(\widehat{\boldsymbol{\sigma}}_h\cdot\boldsymbol{n}, \widehat{\underline{z}}_h\boldsymbol{n})$ on $\dOh$ as
\begin{alignat*}{3}
\widehat{\boldsymbol{\sigma}}_h\cdot\boldsymbol{n}&:=\boldsymbol{\sigma}_h\cdot\boldsymbol{n},
\\
\widehat{\underline{z}}_h\boldsymbol{n}&:=\underline{z}_h\boldsymbol{n},
\end{alignat*}
and impose the transmission conditions
\begin{alignat*}{2}
\langle\widehat{{\boldsymbol{\sigma}}}_h\cdot{\boldsymbol{n}},\mu\rangle_{\dOh\setminus\dO}=&0
&&\quad\forall\;\mu\in M_h,
\\
\langle\widehat{{\boldsymbol{q}}}_h, \boldsymbol{\mu}\rangle_{\dOh\setminus\dO}=&0
&&\quad\forall\;\boldsymbol{\mu} \in \boldsymbol{M}_h.
\end{alignat*}
\end{itemize}

\subsection{Equivalence of the HDG and the WG methods}

\subsubsection{The 2013 WG method \cite{MuWangWangYebiharmonic13} is an HDG method}
We show that the 2013 WG method \cite{MuWangWangYebiharmonic13}  is an HDG method based on the first formulation.
To do that, we first eliminate the unknowns $\boldsymbol{q}_h$
and $\boldsymbol{\sigma}_h$ from the equations defining the HDG methods. We then impose the boundary  and the transmission conditions. After defining the general form of the numerical traces, we finally display the choices of local spaces and the stabilization function which results in the  2014 WG method \cite{MuWangWangYebiharmonic13}. We proceed in several steps.

\subsubsection*{Step 1: The spaces}
We assume that
\[
\boldsymbol{Q}_h=\boldsymbol{\Sigma}_h =:\boldsymbol{V}_h
\quad
\text{ and }
\quad
Z_h=W_h,
\]
and take {\em both} $\mathsf{u}_h:=({u_h},{\uhat_h})$ and $\mathsf{z}_h:=({z_h},{\widehat{z}_h})$ in the space $W_h\times M_h$, where
\[
M_h:=\{\mu\in L^2(\mathcal{F}_h): \; \mu|_F\in M(F)\;\forall F\in \mathcal{F}_h\}.
\]
We set $M_h^0:=\{\mu\in M_h: \mu|_{\dO}=0\}$.

\subsubsection*{Step 2: The auxiliary mapping}
To eliminate the unknowns $\boldsymbol{q}_h$
and $\boldsymbol{\sigma}_h$ from the equations, we express them in terms of an auxiliary linear mapping of 
$\mathsf{u}_h:=({u_h},{\uhat_h})$ and $\mathsf{z}_h:=({z_h},{\widehat{z}_h})$
defined by using the third and first equations, respectively.

Thus, for any $\mathsf{w}:=(\omega,\widehat{\omega})\in L^2(\Oh)\times L^2(\dOh)$, we define ${{\bld{G}_{\mathsf{w}}}}\in \bld{V}_h$ to be the solution of
\begin{alignat*}{2}
({{\bld{G}_{\mathsf{w}}}},{\bld v})_{\Oh} 
&= -(\omega,\nabla \cdot {\bld v})_{\Oh} 
  + \langle {\widehat{\omega}},{\bld v} \cdot {\bld n}\rangle_{\dOh}
&&\quad\forall  {\bld{v}}\in \bld{V}_h.
\end{alignat*}
In this manner, since $\boldsymbol{Q}_h=\boldsymbol{\Sigma}_h=\boldsymbol{V_h}$, when $\mathsf{w}:=\mathsf{u}_h=(u_h,\uhat_h)$, we have that 
$\bld{q}_h={{\bld{G}_{\mathsf{u}_h}}}$, by the third equation, and that, when $\mathsf{w}:=\mathsf{z}_h=(z_h,\widehat{z}_h)$, we have that 
$\bld{\sigma}_h={{\bld{G}_{\mathsf{z}_h}}},$ by the first equation.

\subsubsection*{Step 3: Rewriting the second and fourth equations}
Now, let us work on the remaining equations defining the HDG method. Simple integration by parts in the second and fourth equations gives
\begin{alignat*}{1}
(\nabla\cdot \boldsymbol{\sigma}_h,\omega)_\Oh+\langle (\widehat{{\boldsymbol{\sigma}}}_h-\boldsymbol{\sigma}_h)\cdot{\boldsymbol{n}},
\omega\rangle_\dOh=&(f,\omega)_\Oh,
\\
(\nabla\cdot\boldsymbol{q}_h, s)_{\Oh}+\langle
(\widehat{{\boldsymbol{q}}}_h-\boldsymbol{q}_h)\cdot{\boldsymbol{n}}, s \rangle_{\dOh}=&(z_h,s)_{\Oh},
\end{alignat*}
and, by definition of $\bld{G}_{\mathsf{w}}$ with $\mathsf{w}:=(\omega,\widehat{\omega})$
and $\bld{v}:=\boldsymbol{\sigma}_h$, and
$\bld{G}_{\mathsf{s}}$ with $\mathsf{s}:=(s,\widehat{s}\,)$ and $\bld{v}:=\boldsymbol{q}_h$, we get that
\begin{alignat*}{1}
-(\boldsymbol{G}_{\mathsf{w}}, \boldsymbol{\sigma}_h)_\Oh
+\langle \boldsymbol{\sigma}_h\cdot{\boldsymbol{n}}, \widehat{\omega}\;\rangle_\dOh
+\langle (\widehat{{\boldsymbol{\sigma}}}_h-\boldsymbol{\sigma}_h)\cdot{\boldsymbol{n}},
\omega\rangle_\dOh=&(f,\omega)_\Oh,
\\
-(\boldsymbol{G}_{\mathsf{s}}, \boldsymbol{q}_h)_\Oh
+\langle \boldsymbol{q}_h\cdot{\boldsymbol{n}}, \widehat{s}\;\rangle_\dOh+\langle
(\widehat{{\boldsymbol{q}}}_h-\boldsymbol{q}_h)\cdot{\boldsymbol{n}}, s \rangle_{\dOh}=&(z_h,s)_{\Oh}.
\end{alignat*}
This implies that
\begin{alignat*}{1}
-(\boldsymbol{G}_{\mathsf{w}}, \boldsymbol{\sigma}_h)_\Oh
+\langle (\widehat{{\boldsymbol{\sigma}}}_h-\boldsymbol{\sigma}_h)\cdot{\boldsymbol{n}},
\omega-\widehat{\omega}\;\rangle_\dOh=&(f,\omega)_\Oh-\langle \widehat{\boldsymbol{\sigma}}_h\cdot{\boldsymbol{n}}, \widehat{\omega}\;\rangle_\dOh,
\\
-(\boldsymbol{G}_{\mathsf{s}}, \boldsymbol{q}_h)_\Oh
+\langle
(\widehat{{\boldsymbol{q}}}_h-\boldsymbol{q}_h)\cdot{\boldsymbol{n}}, s-\widehat{s}\;\rangle_{\dOh}=&(z_h,s)_{\Oh}-\langle \widehat{\boldsymbol{q}}_h\cdot{\boldsymbol{n}}, \widehat{s}\;\rangle_\dOh,
\end{alignat*}
and, since  
$\bld{\sigma}_h={{\bld{G}_{\mathsf{z}_h}}}$ and  $\bld{q}_h={{\bld{G}_{\mathsf{u}_h}}}$,  that
\begin{alignat*}{1}
-(\boldsymbol{G}_{\mathsf{w}}, {{\bld{G}_{\mathsf{z}_h}}})_\Oh
+\langle (\widehat{{\boldsymbol{\sigma}}}_h-\boldsymbol{\sigma}_h)\cdot{\boldsymbol{n}},
\omega-\widehat{\omega}\;\rangle_\dOh=&(f,\omega)_\Oh-\langle \widehat{\boldsymbol{\sigma}}_h\cdot{\boldsymbol{n}}, \widehat{\omega}\;\rangle_\dOh,
\\
-(\boldsymbol{G}_{\mathsf{s}}, {{\bld{G}_{\mathsf{u}_h}}})_\Oh
+\langle
(\widehat{{\boldsymbol{q}}}_h-\boldsymbol{q}_h)\cdot{\boldsymbol{n}}, s-\widehat{s}\;\rangle_{\dOh}=&(z_h,s)_{\Oh}-\langle \widehat{\boldsymbol{q}}_h\cdot{\boldsymbol{n}}, \widehat{s}\;\rangle_\dOh,
\end{alignat*}
\subsubsection*{Step 4: Boundary and transmission conditions}
Next, we impose the boundary conditions
\[
\uhat_h=0,
\quad
\widehat{{\boldsymbol{q}}}_h\cdot{\boldsymbol{n}}=0
\qquad\text{ on } \dO,
\]
and the transmissions conditions
\begin{alignat*}{2}
\langle \widehat{{\boldsymbol{\sigma}}}_h\cdot{\boldsymbol{n}}, \widehat{\omega}\;\rangle_{\partial\Oh\setminus\dO}&=0
&&\quad\forall\; \widehat{\omega}\in M_h,
\\
\langle \widehat{{\boldsymbol{q}}}_h\cdot{\boldsymbol{n}}, \widehat{s}\;\rangle_{\partial\Oh\setminus\dO}&=0&&\quad\forall\; \widehat{s}\in M_h,
\end{alignat*}
to get that $\mathsf{u}_h=(u_h,\widehat{u}_h\,)$ is the element of $W_h\times M_h^0$ and
$\mathsf{z}_h=(z_h,\widehat{z_h}\,)$ is the element of $W_h\times M_h$ such that
\begin{alignat*}{1}
-(\boldsymbol{G}_{\mathsf{w}}, {{\bld{G}_{\mathsf{z}_h}}})_\Oh
+\langle (\widehat{{\boldsymbol{\sigma}}}_h-\boldsymbol{\sigma}_h)\cdot{\boldsymbol{n}},
\omega-\widehat{\omega}\;\rangle_\dOh=&(f,\omega)_\Oh,
\\
-(\boldsymbol{G}_{\mathsf{s}}, {{\bld{G}_{\mathsf{u}_h}}})_\Oh
+\langle
(\widehat{{\boldsymbol{q}}}_h-\boldsymbol{q}_h)\cdot{\boldsymbol{n}}, s-\widehat{s}\;\rangle_{\dOh}=&(z_h,s)_{\Oh},
\end{alignat*}
for all $(\omega,\widehat{\omega})\in W_h\times M_h$ and
all $(\omega,\widehat{\omega})\in W_h\times M^0_h$.

\subsubsection*{{\bf Step 5: The numerical traces}}
If we now define the remaining numerical traces by
\begin{alignat*}{2}
\widehat{{\boldsymbol{\sigma}}}_h\cdot{\boldsymbol{n}}
&
:=\boldsymbol{\sigma}_h\cdot{\boldsymbol{n}}
&&\quad\text{ on }\partial\Oh,
\\
\widehat{{\boldsymbol{q}}}_h\cdot{\boldsymbol{n}}&:=\boldsymbol{q}_h\cdot{\boldsymbol{n}}-\tau\, (z_h-\widehat{z}_h)
&&\quad\text{ on }\partial\Oh\setminus\dO,
\end{alignat*}
we get that $\mathsf{u}_h=(u_h,\widehat{u}_h\,)$ is the element of $W_h\times M_h^0$ and
$\mathsf{z}_h=(z_h,\widehat{z_h}\,)$ is the element of $W_h\times M_h$ such that
\begin{alignat*}{2}
-(\boldsymbol{G}_{\mathsf{w}}, \boldsymbol{G}_{\mathsf{z}_h})_\Oh
=&(f,\omega)_\Oh,
\\
-(\boldsymbol{G}_{\mathsf{s}}, \boldsymbol{G}_{\mathsf{u}_h})_\Oh
-\langle
\tau (z_h-\widehat{z}_h), s-\widehat{s}\;\rangle_{\dOh}=&(z_h,s)_{\Oh},
\end{alignat*}
for all $\mathsf{w}=(\omega,\widehat{\omega}\,)\in W_h\times M_h^0$ and
$\mathsf{s}=(s,\widehat{s}\,)\in W_h\times M_h$. 

\subsubsection*{{\bf Step 6: Conclusion}} For the local spaces
\[
\boldsymbol{V}(K):= [\mathcal{P}_k(K)]^d + {\bld x} \mathcal{P}_k(K),
\quad
W(K):=\mathcal{P}_k(K)
\quad\mbox{ and }\quad
M(F):=\mathcal{P}_k(F),
\]
and the stabilization function 
$
\tau|_{\partial K}:=h_K\;\forall\;K\in \Oh,
$ 
this is precisely the 2013 WG method in \cite{MuWangWangYebiharmonic13}; see equations (6) therein.
Thus the
2013 WG method in \cite{MuWangWangYebiharmonic13} is an HDG method.

\subsubsection{The 2014 WG method \cite{MuWangYebiharmonic14} is an HDG method}
We show that the 2014 WG method \cite{MuWangYebiharmonic14}  is an 
HDG method based on the first formulation. To do that, we first eliminate the unknowns 
$\boldsymbol{q}_h$, $\boldsymbol{\sigma}_h$ {\em and} $z_h$ from the equations defining the HDG 
methods. We then impose the boundary  and the transmission conditions. After defining the general form of the numerical traces, we finally display the choices of local spaces and the stabilization function which results in the  2014 WG method \cite{MuWangYebiharmonic14}. We proceed in several steps.

\subsubsection*{Step 1: The spaces}
We assume that
\[
\boldsymbol{Q}_h=\boldsymbol{\Sigma}_h =:\boldsymbol{V}_h,
\]
that
\[
\nabla W(K) \subset \boldsymbol{V}(K):=\boldsymbol{Q}(K)=\boldsymbol{\Sigma}(K)\quad\forall K\in \Oh,
\]
and take the numerical traces 
$(\widehat{\boldsymbol{q}}_h, {\uhat_h})$ in the space $\boldsymbol{N}_h\times M_h$ where
\begin{alignat*}{1}
{\boldsymbol{N}}_h\!&:=\{\boldsymbol{\nu}\in\bld{L}^2(\mathcal{F}_h):\boldsymbol{\nu}\cdot\bld{n}|_{F}\in N(F)\forall F\in \mathcal{F}_h\},
\\
M_h&:=\{\mu\in L^2(\mathcal{F}_h): \; \mu|_F\in M(F)\;\forall F\in \mathcal{F}_h\},
     \end{alignat*}
We set ${\boldsymbol{N}}_h^0:=\{\boldsymbol{\mu}\in {\boldsymbol{N}}_h: \boldsymbol{\mu}\cdot\boldsymbol{n}|_{\dO}=0\}$ and $M_h^0:=\{\mu\in M_h: \mu|_{\dO}=0\}$.

\subsubsection*{Step 2: The auxiliary mappings}
To capture the third equation defining the HDG methods, we proceed as in the previous subsection and introduce the following mapping.
For any $\mathsf{w}:=(\omega,\widehat{\omega})\in L^2(\Oh)\times L^2(\dOh)$, we define 
${{\bld{G}_{\mathsf{w}}}}\in \bld{V}_h$ to be the solution of
\begin{alignat*}{2}
({{\bld{G}_{\mathsf{w}}}},{\bld v})_{\Oh} 
&= -(\omega,\nabla \cdot {\bld v})_{\Oh} 
  + \langle {\widehat{\omega}},{\bld v} \cdot {\bld n}\rangle_{\dOh}
&&\quad\forall  {\bld{v}}\in \bld{V}_h.
\end{alignat*}
In this way,  since $\boldsymbol{Q}_h=\boldsymbol{\Sigma}_h=\boldsymbol{V_h}$, we have that $\boldsymbol{q}_h=\boldsymbol{G}_{\mathsf{u}_h}$, where
$\mathsf{u}_h:=(u_h,\widehat{u}_h)$. 

To capture the fourth equation defining the HDG equations, we introduce an additional mapping. 
For any $\boldsymbol{\mathsf{E}}:=(\boldsymbol{E},\widehat{\boldsymbol{E}})\in [L^2(\Oh)]^d\times[L^2(\dOh)]^d$, we define 
${D_{\boldsymbol{\mathsf{{E}}}}}\in Z_h$ to be the solution of
\begin{alignat*}{2}
({D_{\boldsymbol{\mathsf{{E}}}}},s)_{\Oh} 
&= -(\boldsymbol{E}, \nabla s)_{\Oh} 
  + \langle {\widehat{\boldsymbol{E}}}\cdot\boldsymbol{n}, s\rangle_{\dOh}
&&\quad\forall  s \in Z_h.
\end{alignat*}
In this way, for $\boldsymbol{\mathsf{E}}:=\boldsymbol{\mathsf{Q}_h}:=(\boldsymbol{q}_h, \widehat{\boldsymbol{q}}_h)=(\boldsymbol{G}_{\mathsf{u}_h}, \widehat{\boldsymbol{q}}_h)$,  we have that $z_h=D_{\boldsymbol{\mathsf{Q}_h}}$, by the fourth equations defining the HDG method.

\subsubsection*{Step 3: Rewriting the second equation}
Now, let us work on the second equation defining the HDG method.
By proceeding exactly as in the previous case, we get that
\begin{alignat*}{1}
(f,\omega)_{\Oh}
&=
-(\boldsymbol{G}_{\mathsf{w}}, \boldsymbol{\sigma}_h)_\Oh
+\langle \widehat{\boldsymbol{\sigma}}_h\cdot{\boldsymbol{n}}, \widehat{\omega}\;\rangle_\dOh+\langle (\widehat{{\boldsymbol{\sigma}}}_h-\boldsymbol{\sigma}_h)\cdot{\boldsymbol{n}},
\omega-\widehat{\omega}\;\rangle_\dOh,
\end{alignat*}
with $\mathsf{w}:=(\omega,\widehat{\omega})$. By the first equation with $\boldsymbol{m}:=\boldsymbol{G}_{\mathsf{w}}\in \boldsymbol{Q}_h$, we obtain that
\begin{alignat*}{1}
(f,\omega)_{\Oh}
=&
(z_h, \nabla\cdot\boldsymbol{G}_{\mathsf{w}})_\Oh
-\langle \widehat{z}_h, \boldsymbol{G}_{\mathsf{w}}\cdot{\boldsymbol{n}}\rangle_\dOh
\\
&
+\langle \widehat{\boldsymbol{\sigma}}_h\cdot{\boldsymbol{n}}, \widehat{\omega}\;\rangle_\dOh
+\langle (\widehat{{\boldsymbol{\sigma}}}_h-\boldsymbol{\sigma}_h)\cdot{\boldsymbol{n}},
\omega-\widehat{\omega}\;\rangle_\dOh
\\
=&
-(\nabla z_h, \boldsymbol{G}_{\mathsf{w}})_\Oh
-\langle \widehat{z}_h- z_h, \boldsymbol{G}_{\mathsf{w}}\cdot{\boldsymbol{n}}\rangle_\dOh
\\
&
+\langle \widehat{\boldsymbol{\sigma}}_h\cdot{\boldsymbol{n}}, \widehat{\omega}\;\rangle_\dOh
+\langle (\widehat{{\boldsymbol{\sigma}}}_h-\boldsymbol{\sigma}_h)\cdot{\boldsymbol{n}},
\omega-\widehat{\omega}\;\rangle_\dOh,
\end{alignat*}
and, by the definition of $D_{\boldsymbol{\mathsf{W}}}$ with $\boldsymbol{\mathsf{W}}:=(\boldsymbol{G}_{\mathsf{w}}, \widehat{\boldsymbol{W}})\in \boldsymbol{Q}_h\times\boldsymbol{N}_h$ and $s=z_h$, that
\begin{alignat*}{1}
(f,\omega)_{\Oh}
=&
(D_{\boldsymbol{\mathsf{W}}},  z_h)_\Oh 
-\langle z_h, \widehat{\boldsymbol{W}}\cdot{\boldsymbol{n}}\rangle_\dOh
-\langle \widehat{z}_h-z_h, \boldsymbol{G}_{\mathsf{w}}\cdot{\boldsymbol{n}}\rangle_\dOh
\\
&
+\langle \widehat{\boldsymbol{\sigma}}_h\cdot{\boldsymbol{n}}, \widehat{\omega}\;\rangle_\dOh
+\langle (\widehat{{\boldsymbol{\sigma}}}_h-\boldsymbol{\sigma}_h)\cdot{\boldsymbol{n}},
\omega-\widehat{\omega}\;\rangle_\dOh
\\
=&
(D_{\boldsymbol{\mathsf{W}}}, D_{\boldsymbol{\mathsf{Q}}_h})_\Oh 
\\
&-\langle \widehat{z}_h, \widehat{\boldsymbol{W}}\cdot{\boldsymbol{n}}\rangle_\dOh
-\langle \widehat{z}_h- z_h, (\boldsymbol{G}_{\mathsf{w}}-\widehat{\boldsymbol{W}})\cdot{\boldsymbol{n}}\rangle_\dOh
\\
&
+\langle \widehat{\boldsymbol{\sigma}}_h\cdot{\boldsymbol{n}}, \widehat{\omega}\;\rangle_\dOh
+\langle (\widehat{{\boldsymbol{\sigma}}}_h-\boldsymbol{\sigma}_h)\cdot{\boldsymbol{n}},
\omega-\widehat{\omega}\;\rangle_\dOh,
\end{alignat*}
since $z_h=D_{\boldsymbol{\mathsf{Q}}_h}$.

\subsubsection*{Step 4: Boundary and transmission conditions}
Next,  we impose the boundary conditions
\[
\uhat_h=0,
\quad
\widehat{{\boldsymbol{q}}}_h\cdot{\boldsymbol{n}}=0
\qquad
\text{ on } \dO,
\]
and the transmissions conditions
\begin{alignat*}{2}
\langle \widehat{z}_h, \widehat{\boldsymbol{W}}\cdot{\boldsymbol{n}}\;\rangle_{\partial\Oh\setminus\dO}&=0
&&\quad\forall\; \widehat{\boldsymbol{W}}\in \boldsymbol{N}_h,
\\
\langle \widehat{{\boldsymbol{\sigma}}}_h\cdot{\boldsymbol{n}}, \widehat{\omega}\;\rangle_{\partial\Oh\setminus\dO}&=0
&&\quad\forall\; \widehat{\omega}\in M_h,
\end{alignat*}
to get that $(u_h,\widehat{u}_h,\widehat{\boldsymbol{q}}_h\,) $ is the element of $W_h\times M_h^0
\times  \boldsymbol{N}_h^0$ such that
\begin{alignat*}{3}
(f,\omega)_{\Oh}
=&
( D_{\boldsymbol{\mathsf{W}}},  D_{\boldsymbol{\mathsf{Q}}_h})_\Oh 
&&
-\langle \widehat{z}_h- z_h, (\boldsymbol{G}_{\mathsf{w}}-\widehat{\boldsymbol{W}})\cdot{\boldsymbol{n}}\rangle_\dOh
\\
&&&
+\langle (\widehat{{\boldsymbol{\sigma}}}_h-\boldsymbol{\sigma}_h)\dot{\boldsymbol{n}},
\omega-\widehat{\omega}\;\rangle_\dOh,
\end{alignat*}
for all $(\omega,\widehat{\omega},\widehat{\boldsymbol{W}})\in W_h\times M^0
\times  \boldsymbol{N}_h^0$.

\subsubsection*{{\bf Step 5: The numerical traces}}
It remains to work on the term
\begin{alignat*}{1}
\Theta_h:=
&
-\langle \widehat{z}_h-P z_h, (\boldsymbol{G}_{\mathsf{w}}-\widehat{\boldsymbol{W}})\cdot{\boldsymbol{n}}\rangle_\dOh
+\langle (\widehat{{\boldsymbol{\sigma}}}_h-\boldsymbol{\sigma}_h)\cdot{\boldsymbol{n}},
\omega-\widehat{\omega}\;\rangle_\dOh,
\end{alignat*}
in order to render it a stabilization term by suitably defining the numerical traces $\widehat{z}_h$ and $\widehat{{\boldsymbol{\sigma}}}_h\cdot{\boldsymbol{n}}
$ on $\dOh$.

To do that, we take advantage that we are assuming that 
\[
\nabla W(K)\subset \boldsymbol{Q}(K)\qquad\forall K\in \Oh,
\]
to rewrite the mapping $\mathsf{w}:=(\omega,\widehat{\omega})\mapsto \boldsymbol{G}_{\mathsf{w}}$ as
\[
\boldsymbol{G}_{\mathsf{w}}=\nabla w+ \Phi(\widehat{\omega}-\omega),
\]
where, on each element $K\in \Oh$, 
we define $\Phi(\mu)$ as the element of $\boldsymbol{V}(K)$ such that
\[
(\Phi(\mu), \bld{v})_K=\langle\mu, \bld{v}\cdot\boldsymbol{n}\rangle_{\partial K}
\qquad\forall \bld{v}\in \boldsymbol{V}(K).
\]
Using the definition of $\Phi$, we can now write that
\begin{alignat*}{1}
\Theta_h
:=
&
-\langle \widehat{z}_h- z_h, (\nabla \omega-\widehat{\boldsymbol{W}})\cdot{\boldsymbol{n}}\rangle_\dOh
+\langle (\widehat{{\boldsymbol{\sigma}}}_h-\boldsymbol{\sigma}_h)\cdot{\boldsymbol{n}},
\omega-\widehat{\omega}\;\rangle_\dOh
\\
&
-\langle \widehat{z}_h- z_h, \Phi(\widehat{\omega}-\omega)\cdot{\boldsymbol{n}}\rangle_\dOh
\\
=
&
-\langle \widehat{z}_h- z_h, (\nabla \omega-\widehat{\boldsymbol{W}})\cdot{\boldsymbol{n}}\rangle_\dOh
+\langle (\widehat{{\boldsymbol{\sigma}}}_h-\boldsymbol{\sigma}_h)\cdot{\boldsymbol{n}},
\omega-\widehat{\omega}\;\rangle_\dOh
\\
&
-(\Phi(\widehat{z}_h- z_h), \Phi(\widehat{\omega}-\omega))_{\dOh}
\\
=
&
-\langle \widehat{z}_h- z_h, (\nabla \omega-\widehat{\boldsymbol{W}})\cdot{\boldsymbol{n}}\rangle_\dOh
+\langle (\widehat{{\boldsymbol{\sigma}}}_h-\boldsymbol{\sigma}_h)\cdot{\boldsymbol{n}},
\omega-\widehat{\omega}\;\rangle_\dOh
\\
&
-\langle \Phi(\widehat{z}_h- z_h)\cdot\boldsymbol{n}, \widehat{\omega}-\omega \rangle_\dOh
\\
=&
-\langle \widehat{z}_h- z_h, (\nabla \omega-\widehat{\boldsymbol{W}})\cdot{\boldsymbol{n}}\rangle_\dOh
\\
&
+\langle (\widehat{{\boldsymbol{\sigma}}}_h-\boldsymbol{\sigma}_h)\cdot{\boldsymbol{n}}+\Phi(\widehat{z}_h- z_h)\cdot\boldsymbol{n},
\omega-\widehat{\omega}\;\rangle_\dOh.
\end{alignat*}
This term becomes the stabilization term 
\begin{alignat*}{1}
\Theta_h
=&
\langle \tau_1 (\nabla u_h-\widehat{\boldsymbol{q}}_h)\cdot\boldsymbol{n}, (\nabla \omega-\widehat{\boldsymbol{W}})\cdot{\boldsymbol{n}}\rangle_\dOh
+\langle \tau_2 (u_h-\widehat{u}_h),
\omega-\widehat{\omega}\;\rangle_\dOh.
\end{alignat*}
if we define the numerical traces $(\widehat{z}_n, \widehat{{\boldsymbol{\sigma}}}_h\cdot{\boldsymbol{n}})$ as
\begin{alignat*}{2}
\widehat{z}_h
:=& z_h-\tau_1 (\nabla u_h-\widehat{\boldsymbol{q}}_h)\cdot\boldsymbol{n}
&&\quad\text{ on } \dOh,
\\
\widehat{{\boldsymbol{\sigma}}}_h\cdot{\boldsymbol{n}}
:=&\boldsymbol{\sigma}_h\cdot{\boldsymbol{n}}
-\Phi(\widehat{z}_h- z_h)\cdot\boldsymbol{n}
+\tau_2 (u_h-\widehat{u}_h)
&&\quad\text{ on } \dOh.
\end{alignat*}
This definition certainly allows for the hybridization of the method.

With this choice, we see that the HDG method defines $(u_h,\widehat{u}_h,\widehat{\boldsymbol{q}}_h\,) $ as the element of $W_h\times M_h^0
\times  \boldsymbol{N}_h^0$ such that
\begin{alignat*}{3}
(f,\omega)_{\Oh}
=&
( D_{\boldsymbol{\mathsf{W}}},  D_{\boldsymbol{\mathsf{Q}}_h})_\Oh 
\\&
+\langle \tau_1 (\nabla u_h-\widehat{\boldsymbol{q}}_h)\cdot\boldsymbol{n}, (\nabla \omega-\widehat{\boldsymbol{W}})\cdot{\boldsymbol{n}}\rangle_\dOh
\\
&+\langle \tau_2 (u_h-\widehat{u}_h),
\omega-\widehat{\omega}\;\rangle_\dOh,
\end{alignat*}
for all $(\omega,\widehat{\omega},\widehat{\boldsymbol{W}})\in W_h\times M^0
\times  \boldsymbol{N}_h^0$.

\subsubsection*{{\bf Step 6: Conclusion}}
For the local spaces
\begin{alignat*}{1}
&Z_h:=P_{k-2}(K),
\quad
\bld{V}(K):=[\mathcal{P}_{k-1}(K)]^d,
\quad
W(K):=\mathcal{P}_k(K),
\\
&M(F):=\mathcal{P}_k(F),
\quad
N(F):=\mathcal{P}_{k-1}(F).
\end{alignat*}
and the stabilization functions 
\[
\tau_1|_{\partial K}:=h_K^{-1},
\quad 
\tau_2|_{\partial K}:= h_K^{-3}
\qquad\forall K\in \Oh,
\]
this is precisely the 
2014 WG method \cite{MuWangYebiharmonic14}.
This shows that the 2014 WG method \cite{MuWangYebiharmonic14} is an HDG method.

\subsubsection{The 2014 WG method \cite{WangWangbiharmonic14} is an HDG method}
We show that the 2014 WG method \cite{WangWangbiharmonic14}  is an 
HDG method based on the second formulation. To do that, we proceed exactly as in the previous subsection. So, we first eliminate the unknowns 
$\boldsymbol{q}_h$, $\boldsymbol{\sigma}_h$ {\em and} $\underline{\boldsymbol{z}}_h$ from the equations defining the HDG 
methods. We then impose the boundary  and the transmission conditions. After defining the general form of the numerical traces, we finally display the choices of local spaces and the stabilization function which results in the  2014 WG method \cite{WangWangbiharmonic14}. We proceed in several steps.

\subsubsection*{Step 1: The spaces}
We assume that
\[
\boldsymbol{Q}_h=\boldsymbol{\Sigma}_h =:\boldsymbol{V}_h,
\]
that
\[
\nabla W(K) \subset \boldsymbol{V}(K):=\boldsymbol{Q}(K)=\boldsymbol{\Sigma}(K)\quad\forall K\in \Oh,
\]
and take the numerical traces 
$(\widehat{\boldsymbol{q}}_h, {\uhat_h})$ in the space $\boldsymbol{M}_h\times M_h$ where
\begin{alignat*}{1}
{\boldsymbol{M}}_h\!&:=\{\boldsymbol{\nu}\in\bld{L}^2(\mathcal{F}_h):\boldsymbol{\nu}|_{F}\in \boldsymbol{M}(F)\forall F\in \mathcal{F}_h\},
\\
M_h&:=\{\mu\in L^2(\mathcal{F}_h): \; \mu|_F\in M(F)\;\forall F\in \mathcal{F}_h\},
     \end{alignat*}
We set ${\boldsymbol{M}}_h^0:=\{\boldsymbol{\mu}\in {\boldsymbol{M}}_h: \boldsymbol{\mu}|_{\dO}=\boldsymbol{0}\}$ and $M_h^0:=\{\mu\in M_h: \mu|_{\dO}=0\}$.

\subsubsection*{Step 2: The auxiliary mappings} 
To capture the first and third equations defining the HDG methods, we introduce the following mapping.
For any $\mathsf{w}:=(\omega,\widehat{\omega})\in L^2(\Oh)\times L^2(\dOh)$, we define 
${{\bld{G}_{\mathsf{w}}}}\in \bld{V}_h$ to be the solution of
\begin{alignat*}{2}
({{\bld{G}_{\mathsf{w}}}},{\bld v})_{\Oh} 
&= -(\omega,\nabla \cdot {\bld v})_{\Oh} 
  + \langle {\widehat{\omega}},{\bld v} \cdot {\bld n}\rangle_{\dOh}
&&\quad\forall  {\bld{v}}\in \bld{V}_h.
\end{alignat*}
In this way, we have that $\boldsymbol{q}_h:=\boldsymbol{G}_{\mathsf{u}_h}$, where
$\mathsf{u}_h:=(u_h,\widehat{u}_h)$.

To capture the fourth equation defining the HDG equations, we introduce an additional mapping. 
For any $\boldsymbol{\mathsf{E}}:=(\boldsymbol{E},\widehat{\boldsymbol{E}})\in [L^2(\Oh)]^d\times [L^2(\dOh)]^d$, we define 
${\underline{\boldsymbol{\mathsf{G}}}_{\boldsymbol{\mathsf{{E}}}}}\in W_h$ to be the solution of
\begin{alignat*}{2}
(\underline{\boldsymbol{\mathsf{G}}}_{\boldsymbol{\mathsf{{E}}}},w)_{\Oh} 
&= -(\boldsymbol{E}, \nabla\cdot\underline{\boldsymbol{s}})_{\Oh} 
  + \langle {\widehat{\boldsymbol{E}}}, \underline{\boldsymbol{s}}\boldsymbol{n}\rangle_{\dOh}
&&\quad\forall  \underline{\boldsymbol{s}} \in \underline{\boldsymbol{Z}}_h.
\end{alignat*}
In this way, for $\boldsymbol{\mathsf{E}}:=\boldsymbol{\mathsf{Q}_h}:=(\boldsymbol{q}_h, \widehat{\boldsymbol{q}}_h)=(\boldsymbol{G}_{\mathsf{u}_h},\widehat{\boldsymbol{q}}_h)$,  we have that $\underline{\boldsymbol{z}}_h=\underline{\boldsymbol{\mathsf{G}}}_{\boldsymbol{\mathsf{Q}_h}}$, by the fourth equation defining the HDG method.

\subsubsection*{Step 3: Rewriting the second equation}
Now, let us work on the second equation defining the HDG method. Proceeding exactly and in the two previous cases, we get
\begin{alignat*}{1}
(f,\omega)_{\Oh}
&=
-(\boldsymbol{G}_{\mathsf{w}}, \boldsymbol{\sigma}_h)_\Oh
+\langle \widehat{\boldsymbol{\sigma}}_h\cdot{\boldsymbol{n}}, \widehat{\omega}\;\rangle_\dOh+\langle (\widehat{{\boldsymbol{\sigma}}}_h-\boldsymbol{\sigma}_h)\cdot{\boldsymbol{n}},
\omega-\widehat{\omega}\;\rangle_\dOh,
\end{alignat*}
with $\mathsf{w}:=(\omega,\widehat{\omega})$. By the first equation with $\boldsymbol{m}:=\boldsymbol{G}_{\mathsf{w}}$, we obtain that
\begin{alignat*}{1}
(f,\omega)_{\Oh}
=&
(\underline{\boldsymbol{z}}_h, \nabla \boldsymbol{G}_{\mathsf{w}})_\Oh
-\langle \widehat{\underline{\boldsymbol{z}}}_h{\boldsymbol{n}}, \boldsymbol{G}_{\mathsf{w}}\rangle_\dOh
\\
&
+\langle \widehat{\boldsymbol{\sigma}}_h\cdot{\boldsymbol{n}}, \widehat{\omega}\;\rangle_\dOh
+\langle (\widehat{{\boldsymbol{\sigma}}}_h-\boldsymbol{\sigma}_h)\cdot{\boldsymbol{n}},
\omega-\widehat{\omega}\;\rangle_\dOh
\\
=&
-(\nabla\cdot  \underline{\boldsymbol{z}}_h, \boldsymbol{G}_{\mathsf{w}})_\Oh
-\langle (\widehat{\underline{\boldsymbol{z}}}_h-\underline{\boldsymbol{z}}_h){\boldsymbol{n}}, \boldsymbol{G}_{\mathsf{w}}\rangle_\dOh
\\
&
+\langle \widehat{\boldsymbol{\sigma}}_h\cdot{\boldsymbol{n}}, \widehat{\omega}\;\rangle_\dOh
+\langle (\widehat{{\boldsymbol{\sigma}}}_h-\boldsymbol{\sigma}_h)\cdot{\boldsymbol{n}},
\omega-\widehat{\omega}\;\rangle_\dOh,
\end{alignat*}
and, by the definition of $\underline{\boldsymbol{\mathsf{G}}}_{\boldsymbol{\mathsf{{W}}}}$
with $\boldsymbol{\mathsf{W}}:=(\boldsymbol{G}_{\mathsf{w}}, \widehat{\boldsymbol{W}})$
and $\underline{\boldsymbol{s}}:=\underline{\boldsymbol{z}}_h$, that
\begin{alignat*}{1}
(f,\omega)_{\Oh}
=&
(\underline{\boldsymbol{\mathsf{G}}}_{\boldsymbol{\mathsf{{W}}}},   \underline{\boldsymbol{z}}_h)_\Oh
-\langle \underline{\boldsymbol{z}}_h{\boldsymbol{n}}, \widehat{\boldsymbol{W}}\rangle_\dOh
-\langle (\widehat{\underline{\boldsymbol{z}}}_h- \underline{\boldsymbol{z}}_h){\boldsymbol{n}}, \boldsymbol{G}_{\mathsf{w}}\rangle_\dOh
\\
&
+\langle \widehat{\boldsymbol{\sigma}}_h\cdot{\boldsymbol{n}}, \widehat{\omega}\;\rangle_\dOh
+\langle (\widehat{{\boldsymbol{\sigma}}}_h-\boldsymbol{\sigma}_h)\cdot{\boldsymbol{n}},
\omega-\widehat{\omega}\;\rangle_\dOh
\\
=&
( \underline{\boldsymbol{\mathsf{G}}}_{\boldsymbol{\mathsf{{W}}}},  \underline{\boldsymbol{\mathsf{G}}}_{\boldsymbol{\mathsf{{Q}}}_h})_\Oh
\\
&-\langle \widehat{\underline{\boldsymbol{z}}}_h{\boldsymbol{n}}, \widehat{\boldsymbol{W}}\rangle_\dOh
-\langle (\widehat{\underline{\boldsymbol{z}}}_h- \underline{\boldsymbol{z}}_h){\boldsymbol{n}}, \boldsymbol{G}_{\mathsf{w}}-\widehat{\boldsymbol{W}}\rangle_\dOh
\\
&
+\langle \widehat{\boldsymbol{\sigma}}_h\cdot{\boldsymbol{n}}, \widehat{\omega}\;\rangle_\dOh
+\langle (\widehat{{\boldsymbol{\sigma}}}_h-\boldsymbol{\sigma}_h)\cdot{\boldsymbol{n}},
\omega-\widehat{\omega}\;\rangle_\dOh
\end{alignat*}
since $\underline{\boldsymbol{z}}_h=\underline{\boldsymbol{\mathsf{G}}}_{\boldsymbol{\mathsf{{Q}}}_h}$.

\subsubsection*{Step 4: Boundary and transmission conditions}
Next,  we impose the boundary conditions
\[
\uhat_h=0,
\quad
\widehat{{\boldsymbol{q}}}_h=\boldsymbol{0}\qquad
\text{ on }\dO,
\]
and the transmissions conditions
\begin{alignat*}{2}
\langle \widehat{\underline{\boldsymbol{z}}}_h\boldsymbol{n}, \widehat{\boldsymbol{W}}\;\rangle_{\partial\Oh\setminus\dO}&=0
&&\quad\forall\; \widehat{\boldsymbol{W}}\in \boldsymbol{M}_h,
\\
\langle \widehat{{\boldsymbol{\sigma}}}_h\cdot{\boldsymbol{n}}, \widehat{\omega}\;\rangle_{\partial\Oh\setminus\dO}&=0
&&\quad\forall\; \widehat{\omega}\in M_h,
\end{alignat*}
to get that $(u_h,\widehat{u}_h,\widehat{\boldsymbol{q}})\in W_h\times M_h^0\times
\boldsymbol{M}^0_h$ satisfies
\begin{alignat*}{3}
(f,\omega)_{\Oh}
=&
( \underline{\boldsymbol{\mathsf{G}}}_{\boldsymbol{\mathsf{{W}}}},  \underline{\boldsymbol{\mathsf{G}}}_{\boldsymbol{\mathsf{{Q}}}_h})_\Oh
&&
-\langle (\widehat{\underline{\boldsymbol{z}}}_h- \underline{\boldsymbol{z}}_h){\boldsymbol{n}}, \boldsymbol{G}_{\mathsf{w}}-\widehat{\boldsymbol{W}}\rangle_\dOh
\\
&&&
+\langle (\widehat{{\boldsymbol{\sigma}}}_h-\boldsymbol{\sigma}_h)\cdot{\boldsymbol{n}},
\omega-\widehat{\omega}\;\rangle_\dOh,
\end{alignat*}
for all $(\omega,\widehat{\omega},\widehat{\boldsymbol{W}})\in W_h\times M_h^0\times
\boldsymbol{M}^0_h$.

\subsubsection*{{\bf Step 5: The numerical traces}}
Let us now work on the term
\begin{alignat*}{1}
\Theta_h:=
&
-\langle (\widehat{\underline{\boldsymbol{z}}}_h- \underline{\boldsymbol{z}}_h){\boldsymbol{n}}, \boldsymbol{G}_{\mathsf{w}}-\widehat{\boldsymbol{W}}\rangle_\dOh
+\langle (\widehat{{\boldsymbol{\sigma}}}_h-\boldsymbol{\sigma}_h)\cdot{\boldsymbol{n}},
\omega-\widehat{\omega}\;\rangle_\dOh,
\end{alignat*}
in order to render it a stabilization term by suitably defining the numerical traces $\widehat{\underline{\boldsymbol{z}}}_h$ and $\widehat{{\boldsymbol{\sigma}}}_h\cdot{\boldsymbol{n}}
$ on $\dOh$.
To do that, we take advantage that we are assuming that 
\[
\nabla W(K)\subset \boldsymbol{V}(K)\qquad\forall K\in \Oh,
\]
to rewrite the mapping $\mathsf{w}:=(\omega,\widehat{\omega})\mapsto \boldsymbol{G}_{\mathsf{w}}$ as
\[
\boldsymbol{G}_{\mathsf{w}}=\nabla w+ \Phi(\widehat{\omega}-\omega),
\]
where, on the element $K\in \Oh$,  $\Phi(\mu)$ is the element of $\boldsymbol{V}(K)$ such that
\[
(\Phi(\mu), \bld{v})_K=\langle\mu, \bld{v}\cdot\boldsymbol{n}\rangle_{\partial K}
\qquad\forall \bld{v}\in \boldsymbol{V}(K).
\]
We also define, on the element $K\in \Oh$, $\Psi(\boldsymbol{\mu})$ as the element of 
\[
R(\partial K):=\{\mu\in L^2(\partial K):\; \mu|_F\in M(F) + W(K)|_F\;\text{ for all faces $F$ of $K$}\},
\]
defined by
\[
\langle \Psi(\boldsymbol{\mu}), r \rangle_{\partial K}
=
\langle \boldsymbol{\mu}, \Phi(r) \rangle_{\partial K}\qquad \forall r\in R(\partial K).
\]

With the help of the lifting $\Phi$ and the mapping $\Psi$, we can write
\begin{alignat*}{1}
\Theta_h
:=
&
-\langle (\widehat{\underline{\boldsymbol{z}}}_h- \underline{\boldsymbol{z}}_h){\boldsymbol{n}}, \nabla\omega-\widehat{\boldsymbol{W}}\rangle_\dOh
+\langle (\widehat{{\boldsymbol{\sigma}}}_h-\boldsymbol{\sigma}_h)\cdot{\boldsymbol{n}},
\omega-\widehat{\omega}\;\rangle_\dOh
\\
&-\langle (\widehat{\underline{\boldsymbol{z}}}_h- \underline{\boldsymbol{z}}_h){\boldsymbol{n}}, \Phi(\widehat{\omega}-\omega)\rangle_\dOh
\\
=
&
-\langle (\widehat{\underline{\boldsymbol{z}}}_h- \underline{\boldsymbol{z}}_h){\boldsymbol{n}}, \nabla\omega-\widehat{\boldsymbol{W}}\rangle_\dOh
+\langle (\widehat{{\boldsymbol{\sigma}}}_h-\boldsymbol{\sigma}_h)\cdot{\boldsymbol{n}},
\omega-\widehat{\omega}\;\rangle_\dOh
\\
&-\langle \Psi((\widehat{\underline{\boldsymbol{z}}}_h- \underline{\boldsymbol{z}}_h){\boldsymbol{n}}), \widehat{\omega}-\omega\rangle_\dOh.
\\
=
&
-\langle (\widehat{\underline{\boldsymbol{z}}}_h- \underline{\boldsymbol{z}}_h){\boldsymbol{n}}, \nabla\omega-\widehat{\boldsymbol{W}}\rangle_\dOh
\\
&
+\langle (\widehat{{\boldsymbol{\sigma}}}_h-\boldsymbol{\sigma}_h)\cdot{\boldsymbol{n}}+\Psi((\widehat{\underline{\boldsymbol{z}}}_h- \underline{\boldsymbol{z}}_h){\boldsymbol{n}}),
\omega-\widehat{\omega}\;\rangle_\dOh.
\end{alignat*}
This terms becomes the stabilization term
\begin{alignat*}{1}
\Theta_h:=
&
\langle \tau_1(\nabla u_h-\widehat{\boldsymbol{q}}_h), (\nabla \omega-\widehat{\boldsymbol{W}}\rangle_\dOh
+\langle \tau_2(u_h-\widehat{u}_h),
(\omega-\widehat{\omega})\;\rangle_\dOh,
\end{alignat*}
if we define the numerical traces on $\partial\Oh$ by
\begin{alignat*}{2}
\widehat{\underline{\boldsymbol{z}}}_h
:=& \underline{\boldsymbol{z}}_h+\tau_1(\nabla u_h-\widehat{\boldsymbol{q}}_h)
\\
\widehat{{\boldsymbol{\sigma}}}_h\cdot{\boldsymbol{n}}
:=&\boldsymbol{\sigma}_h\cdot{\boldsymbol{n}}
-\Psi((\widehat{\underline{\boldsymbol{z}}}_h-\underline{\boldsymbol{z}}_h)\boldsymbol{n})
+\tau_2(u_h-\widehat{u}_h).
\end{alignat*}
Again, this definition allows for the hybridization of the method.

With this choice, the HDG method defines $(u_h,\widehat{u}_h,\widehat{\boldsymbol{q}})\in W_h\times M_h^0\times
\boldsymbol{M}^0_h$ as the solution of
\begin{alignat*}{3}
(f,\omega)_{\Oh}
=&
( \underline{\boldsymbol{\mathsf{G}}}_{\boldsymbol{\mathsf{{W}}}},  \underline{\boldsymbol{\mathsf{G}}}_{\boldsymbol{\mathsf{{Q}}}_h})_\Oh
&&
+\langle \tau_1(\nabla u_h-\widehat{\boldsymbol{q}}_h), (\nabla \omega-\widehat{\boldsymbol{W}}\rangle_\dOh
\\
&&&
+\langle \tau_2(u_h-\widehat{u}_h),
(\omega-\widehat{\omega})\;\rangle_\dOh,
\end{alignat*}
for all $(\omega,\widehat{\omega},\widehat{\boldsymbol{W}})\in W_h\times M_h^0\times
\boldsymbol{M}^0_h$.

\subsubsection*{{\bf Step 6: Conclusion}}
For the local spaces
\begin{alignat*}{1}
&\underline{\bld{Z}}(K):=[\mathcal{P}_{k-2}(K)]^{d\times d},
\bld{\Sigma}(K)= \bld{Q}(K):=[\mathcal{P}_{k-1}(K)]^d,
\;\;
W(K):=\mathcal{P}_k(K),
\\
&
M(F):=\mathcal{P}_{k-2}(F),
\;\;
\boldsymbol{M}(F):=[\mathcal{P}_{k-2}(F)]^d.
\end{alignat*}
and  the Lehrenfeld-Sch\"oberl  stabilization functions \cite{Lehrenfeld10,Schoeberl09}
\[
\tau_1|_{\partial K}:=h_K^{-1} \boldsymbol{P}_{\boldsymbol{M}} 
\quad
\tau_2|_{\partial K}:=h_K^{-3} P_M
 \quad\text{ on the face $F$ of the element $K\in \Oh$,}
\]
where $h_K$ is the diameter of rhte element $K\in \Oh$, and $P_M$ is the $L^2$-projection into $M_h$ and $\boldsymbol{P}_{\boldsymbol{M}}$ is the $L^2$-projection into $\boldsymbol{M}_h$, this is precisely the 2014 WG method in \cite{WangWangbiharmonic14}. This shows that the 2014 WG method \cite{WangWangbiharmonic14} is an HDG method.

\section{A critique of questionable claims}

\subsection{On the novelty of the discrete weak gradient}
WG authors assert that their discovery of the {\em discrete weak
  gradient} is the main novelty of their work; see, for example, 
lines -9 to -6 in page 765 of \cite{MR3077972}. 

Unfortunately, this is not a novelty. In fact, this very concept was
introduced, under a different name, at least back in 1977: see page 296, 
Remark 2, equation (2.15) in \cite{RaviartThomas77}.
Moreover, in the
unifying work on DG methods for second-order elliptic equations
\cite{ArnoldBrezziCockburnMarini02}, the discrete weak gradients ``invented'' by
the WG authors, were systematically used, again, under a different name, 
to describe all the then-known DG
methods for second-order elliptic equations.
Finally, these ``discrete weak gradients"
(an other similarly defined operators) have been the bread 
and butter  of all DG approximations since the very inception of the method      
back in 1973, see \cite{CockburnShuReview01,CockburnZAMM03,CockburnEncy04}. 
None of these references have ever been quoted in WG papers. 

Other researchers, have also proposed
interpretations of the``discrete weak gradient" similar to the one proposed by
the WG authors,  but for much more fruitful ends. 
In particular, Adrian Lew and
collaborators \cite{EyckLew06}, called the ``discrete weak gradient'', the ``DG
derivative'' and used it to actually define {\em new} DG methods for nonlinear
elasticity. No mention to their work has been
ever made in any WG  paper.

\subsection{On the flexibility for handling elements of arbitrary shape}
WG authors claim that a novelty of their method is its
ability to handle elements of general geometry; see, for example, the last line
of page 5 in \cite{MuWangYeZhao14}. Unfortunately,
this is a direct consequence of the WG method being nothing but a HDG method and hence a DG method; 
it cannot be presented as a novelty.

\subsection{On the Lehrenfeld-Sch\"oberl stabilization}
In the papers \cite{WangYemixedsecondorder14,MuWangYeLS15,WangYeStokes16}, WG authors arrogated for themselves the introduction of the stabilization function
\[
\tau:= \frac{\rho}{h_K}P_{M(F)},
\]
where $\rho$ is a positive constant,  
$h_K$ is the diameter of the element $K$ and $P_{M(F)}$ is the
$L^2$-projection into $M(F)$. As we have pointed out above, this stabilization function
was introduced  by
 Christoph Lehrenfeld back in May of 2010 in his Diploma Thesis entitled  ``Hybrid
Discontinuous Galerkin methods for solving incompressible flow problems",
Aachen, Germany, under the direction of Joachim Sch\"oberl; see Remark 1.2.4 in page 17 in \cite{Lehrenfeld10} and \cite{Schoeberl09}. In his 2015 paper, Oikawa \cite{Oikawa15} does refer to the projection under consideration as the Lehrenfeld-Sch\"oberl stabilization,
showing that this is a well-known fact in the community.

Joachim Sch\"oberl spoke about this stabilization function in his plenary talk at ICOSAHOM '09 conference, June 22-26, Trondheim, Norway.
More importantly, he  spoke about this method  in the Fall 2010 Finite Element Circus 
held at the Institute for Mathematics and its Applications in Minneapolis, MN on November 5-6, 2010 to honor Profs. R. Falk, J. Pasciak and L. Wahlbin. 
As we can see in the list of participants

\

https://www.ima.umn.edu/2010-2011/SW11.5-6.10/\#participants

\

\noindent The authors of the first paper on WG methods  were present. They
were aware of the Lehrenfeld-Sch\"oberl stabilization function and the HDG method by  Christoph Lehrenfeld 
and Joachim Sch\"oberl {\em before} they wrote the papers in which they used it  \cite{WangYemixedsecondorder14,MuWangYeLS15,WangYeStokes16}.

\subsection{On the definition of the Weak Galerkin methods} The definition of what constitutes a WG method has changed with time.

The very first WG methods \cite{WangYeMain13} were devised to discretize the steady-state
second-order elliptic equation
\[
-\nabla\cdot (a\nabla u)+ \nabla\cdot(bu)+ c u=f.
\]
In that paper, the authors write that 

\

{\em The weak Galerkin finite element method, as detailed in Section 4, is
  closely related to the mixed finite element method (see [9-12,14,16]) with
a hybridized interpretation of Fraejis de Veubeque [21,22] {\bf and the
  hybridizable discontinuous Galerkin (HDG) method [18]}. The hybridized
formulation introduces a new term, known as the Lagrange multiplier, on the
boundary of each element. The Lagrange multiplier is known to approximate the
original function $u=u(x)$ on the boundary of each element. {\bf In fact, the WG
schemes introduced in Section 4 are equivalent to some standard mixed methods
and the HDG method when $b=0, c=0$ and $a=const$ in (1.1). But the WG method is
different from these methods when the coefficients are general variable functions.} The concept of {\em
  weak gradients} shall provide a systematic framework for dealing with
discontinuous functions defined on elements and their boundaries in a near
classical sense.} 

\

Here, we have used boldface to mark the sentences added from the original version \cite{WangYeMain11v1}. We can see that in the original version, 
no mention was made to a possible relation between the WG and the HDG (and mixed) methods. This is because in this setting, the WG methods are 
in fact hybridized mixed methods. In the second version of the paper, the authors  claim, incorrectly, that the WG method is not an HDG method when the coefficients are general variable functions, as we have argued in \S 2. 

This gives the impression that the variability of the coefficients is what makes the difference between the WG and the HDG methods. And yet,
when dealing with the Stokes problem and the bilaplacian, whose differential operators present coefficients which are all constant, the WG authors still claim that their WG method is not an HDG method.

Finally, let us point out that for this convection-diffusion-reaction equation, no stabilization is needed due to the presence of the reaction term.  
However, when this term is not present,  the introduction of the stabilization function, typical of HDG methods, is necessary. This is why the WG authors  
incorporated the stabilization function in their subsequent papers.

\section{Conclusion}
We have established that the WG methods are nothing but rewritings of the HDG methods. We have done that 
for steady-state diffusion, (briefly for) the Stokes equations of incompressible flow, and for the bilaplacian equation.
Similar results can be established for {\bf all} other applications. This paper shows that the WG methods are a simple rebranding of the HDG methods.

\

\providecommand{\bysame}{\leavevmode\hbox to3em{\hrulefill}\thinspace}
\providecommand{\MR}{\relax\ifhmode\unskip\space\fi MR }
\providecommand{\MRhref}[2]{%
  \href{http://www.ams.org/mathscinet-getitem?mr=#1}{#2}
}
\providecommand{\href}[2]{#2}

\end{document}